\documentclass[11pt]{amsart}

\usepackage{amsmath}
\usepackage{amssymb}
\usepackage[active]{srcltx}

\newcommand{\R}{{\mathbb R}}       % Field of real numbers
\newcommand{\Z}{{\mathbb Z}}       % Ring of integer numbers
\newcommand{\DD}{{\mathcal D}}

\newcommand{\AD}{{\mathcal {AD}}}
\newcommand{\AS}{{\mathcal {AS}}}
\newcommand{\ST}{{\mathcal {ST}}}

\newcommand{\CC}{{\mathcal C}}

\newcommand{\diam}{{\rm diam}}
\newcommand{\dist}{{\rm dist}}
\newcommand{\ds}{\displaystyle }
\newcommand{\ts}{\textstyle }

\newcommand{\rf}[1]{{(\ref{#1})}}
\newcommand{\sk}{{{\mathcal S}_{(k,i)}}}
\newcommand{\skd}{{{\mathcal S}^*_{(k,i)}}}
\newcommand{\skk}{{\overline{{\mathcal S}}_{(k,i)}}}
\newcommand{\skkd}{{\overline{{\mathcal S}}^*_{(k,i)}}}

\newcommand{\supp}{{\rm supp}}

\newcommand{\vphi}{{\varphi}}
\newcommand{\ve}{{\varepsilon}}

\newcommand{\wt}[1]{{\widetilde{#1}}}
\newcommand{\wh}[1]{{\widehat{#1}}}

\newcommand{\bmo}{{B\!M\!O}}
\newcommand{\rbmo}{{RBMO}}

\newtheorem{theorem}{Theorem}[section]
\newtheorem{lemma}[theorem]{Lemma}

\newtheorem{propo}[theorem]{Proposition}

\theoremstyle{definition}
\newtheorem{definition}[theorem]{Definition}
\newtheorem{example}[theorem]{Example}

\theoremstyle{remark}
\newtheorem{rem}[theorem]{Remark}

\numberwithin{equation}{section}

\newcommand{\brem}{\begin{rem}}
\newcommand{\erem}{\end{rem}}

\newcommand{\bexam}{\begin{example}}
\newcommand{\eexam}{\end{example}}

\begin{document}

\title[Weights for Calder\'on-Zygmund operators]
{Weighted norm inequalities for Calder\'on-Zygmund operators
without doubling conditions}

\author[XAVIER TOLSA]{Xavier Tolsa}

\address{Instituci\'o Catalana de Recerca i Estudis Avan\c{c}ats (ICREA) and Departament de Matem\`atiques,
Universitat Aut\`onoma de Bar\-ce\-lo\-na, Spain}

\email{xtolsa@mat.uab.cat}

\thanks{Partially supported by grants MTM2004-00519 and
Acci\'on Integrada HF2004-0208 (Spain), and 2001-SGR-00431 (Generalitat
de Catalunya)}


\date{}

%\keywords{Calder\'on-Zygmund operators, weights, non
%doubling measures}

\begin{abstract}
Let $\mu$ be a Borel measure on $\R^d$ which may be non doubling.
The only condition that $\mu$ must satisfy is $\mu(B(x,r))\leq Cr^n$
for all $x\in\R^d$, $r>0$ and for some fixed $n$ with $0<n\leq d$.
In this paper we introduce a maximal operator $N$, which coincides
with the maximal Hardy-Littlewood operator if $\mu(B(x,r))\approx r^n$ for
$x\in\supp(\mu)$,
and we show that all $n$-dimensional Calder\'on-Zygmund operators are
bounded on $L^p(w\,d\mu)$ if and only if $N$ is bounded on $L^p(w\,d\mu)$,
for a fixed $p\in(1,\infty)$. Also, we prove that this happens if and
only if some conditions of Sawyer type hold.
We obtain analogous results about the weak $(p,p)$ estimates. This
type of weights do not satisfy a reverse H\"older inequality, in general,
but some kind of self improving property still holds. On the other hand, if $f
\in\rbmo(\mu)$ and $\ve>0$ is small enough,
then $e^{\ve f}$ belongs to this class of weights.
\end{abstract}

\maketitle

\section{Introduction}

Let $\mu$ be some Borel measure on $\R^d$ satisfying
\begin{equation} \label{creixlin}
\mu(B(x,r))\leq C_0\,r^n \qquad \mbox{for all $x\in\R^d$, $r>0$},
\end{equation}
where $n$ is some fixed constant (which may be non integer) with
$0<n\leq d$.
In this paper we obtain a characterization of all the weights $w$ such
that, for every $n$-dimensional Calder\'on-Zygmund operator (CZO) $T$
which is bounded on $L^2(\mu)$, the following
weighted inequality holds:
\begin{equation} \label{int1}
\int |Tf|^p\,w\,d\mu \leq C \int |f|^p\,w\,d\mu\qquad \mbox{for all
$f\in L^p(w)$},
\end{equation}
where $C$ is independent of $f$, $1<p<\infty$, and
$L^p(w):=L^p(w\,d\mu)$.
It is shown that these weights $w$ are those such that a suitable
maximal operator $N$ (defined below) is bounded on $L^p(w)$. We also
prove an analogous result for the weak $(p,p)$ estimates.

Moreover, we show that the $L^p$ weights for CZO's (and for $N$) satisfy
a self improving property. Loosely speaking, weak weighted inequalities
for $w$ and for the dual weight $w^{-1/(p-1)}$ imply strong weighted
inequalities for $w$ and its dual weight.
Let us remark that we do {\em not} assume that the underlying measure
$\mu$ is doubling. Recall that $\mu$ is said to be doubling if there
exists some constant $C$ such that
%$$\mu(B(x,2r)) \leq C\, \mu(B(x,r))\qquad \mbox{for all $x\in\supp(\mu)$,
%$r>0$}.$$
$\mu(B(x,2r)) \leq C\, \mu(B(x,r))$ for all $x\in\supp(\mu)$ and $r>0$.

In the particular case where $\mu$  coincides with the Lebesgue
measure on $\R^d$, it is known that the weighted inequality
\rf{int1} holds for all $d$-dimensional CZO's if
and only if $w$ is an $A_p$ weight. This result was obtained by
Coifman and Fefferman \cite{CF}, and it generalizes a previous result by
Hunt, Muckenhoupt and Wheeden \cite{HMW} about the Hilbert
transform. Let us recall that Muckenhoupt proved \cite{Muck} that the
$A_p$ weights are precisely those weights $w$ for which the
Hardy-Littlewood operator is bounded on $L^p(w)$ (always assuming
$\mu$ to be the Lebesgue measure on $\R^d$). So
the $L^p$ weights for CZO's and the $L^p$ weights for
the maximal Hardy-Littlewood operator coincide in this case.

Suppose now that the measure $\mu$ satisfies
\begin{equation} \label{AD}
\mu(B(x,r))\approx r^n \qquad \mbox{for all $x\in\supp(\mu)$, $r>0$},
\end{equation}
where $A\approx B$ means that there is some constant $C>0$ such that
$C^{-1}A \leq B \leq CA$, with $C$
depending only on $n$ and $d$ (and also on $C_0$ sometimes), in general. In
this case the results (and their proofs) are analogous to the ones for the
Lebesgue measure. Namely, \rf{int1} holds for all $n$-dimensional
CZO's if and only if $w\in A_p$, which is equivalent to say that the
maximal Hardy-Littlewood operator is bounded on $L^p(w)$.

Many other results about weights for CZO's can be found in the
literature. In most of them it is assumed that $\mu$ is either the
Lebesgue measure on $\R^d$ or the underlying measure of a
space of homogeneous type, satisfying \rf{AD}. See for example
\cite{Perez} and the recent result on the two weight problem for the Hilbert transform in \cite{Volberg}.

It is much more difficult to find results where \rf{AD} does not hold.
%In this situation the usual techniques don't work, in general.
Saksman \cite{Saksman} has obtained some results
concerning the weights for the Hilbert transform $H$ on arbitrary
bounded subsets of $\R$
(with $\mu$ being the Lebesgue measure restricted to these subsets).
These results relate the boundedness of $H$ on
$L^p(w)$ with some operator properties of $H$, and
quite often his arguments are of complex analytic nature.
Orobitg and P\'erez \cite{OP} have studied the $A_p$ classes of weights
with respect to arbitrary measures on $\R^d$, which may be non
doubling. In particular, they have shown that if $w$ is an $A_p$
weight, then the centered maximal Hardy-Littlewood operator is bounded
on $L^p(w)$, and that if moreover $\mu$ satisfies \rf{creixlin}, then
all $n$-dimensional CZO's are also bounded on $L^p(w)$. Other more recent result which involve
the operator
$$M_kf(x)=\sup_{x\in B}\frac{1}{\mu(kB)}\int_B |f(y)|d\mu(y)\quad\mbox{for $x\in {\rm supp}(\mu)$ and $k>1$},$$ where the supremum is taken over all balls $B$ containing $x$,  have been obtained in \cite{Komori}.

Our approach uses real variable techniques and it is based on the ideas
and methods developed in \cite{Tolsa2}, \cite{Tolsa3} and
\cite{Tolsa4} to extend Calder\'on-Zygmund theory to the the setting of
non doubling measures. Indeed, recently it has been shown that
the doubling assumption is not essential for many results of
Calder\'on-Zygmund theory. See \cite{NTV1}, \cite{Tolsa1},
\cite{NTV2}, \cite{MMNO} and \cite{GM}, for instance,
in addition to the references cited above.

In order to state our results more precisely, we need to introduce some
definitions.
A kernel $k(\cdot,\cdot):\R^d\times\R^d\to\R$
is called a ($n$-dimensional) Calder\'on-Zyg\-mund (CZ) kernel if
\begin{itemize}
\item[(1)] $\ds |k(x,y)|\leq \frac{C_1}{|x-y|^n}\quad$ if $x\neq y$,
\item[(2)] there exists some fixed constant $0<\gamma\leq1$ such that
$$|k(x,y)-k(x',y)| + |k(y,x)-k(y,x')|
\leq C_2\,\frac{|x-x'|^\gamma}{|x-y|^{n+\gamma}}$$
if $|x-x'|\leq |x-y|/2$.
\end{itemize}
Throughout all the paper we will
assume that $\mu$ is a Radon measure on $\R^d$ satisfying \rf{creixlin}.
We say that $T$ is a ($n$-dimensional) CZO
associated to the kernel
$k(x,y)$ if for any compactly supported function $f\in L^2(\mu)$
\begin{equation} \label{eq***}
Tf(x) = \int k(x,y)\,f(y)\,d\mu(y)\qquad \mbox{if $x\not\in\supp(\mu)$}.
\end{equation}
and $T$ is bounded on $L^2(\mu)$ (see the paragraph below regarding
this question). If we want to make explicit the constant $\gamma$
which appears in the
second property of the CZ kernel, we will write $T\in CZO(\gamma)$.

The integral in \rf{eq***} may be non convergent for $x\in\supp(\mu)$, even
for ``very nice'' functions, such as $\CC^\infty$ functions with compact
support. For this reason it is convenient to introduce the truncated operators
$T_\ve$, $\ve>0$:
$$T_\varepsilon f(x) = \int_{|x-y|>\varepsilon} k(x,y)\,f(y)\,d\mu(y).$$
Then we say that $T$ is bounded on $L^2(\mu)$ if the operators $T_\ve$
are bounded on $L^2(\mu)$ uniformly on $\ve>0$.

Now we will define the maximal operator $N$.
For $0<r< R$ and a fixed $x\in\supp(\mu)$, we consider the
function
$$\vphi_{x,r,R}(y) = \left\{ \begin{array}{cl}
1/r^n & \mbox{if $0\leq |x-y|\leq r$,} \\
1/|x-y|^n & \mbox{if $r\leq |x-y|\leq R$,} \\
0 &  \mbox{if $|x-y|> R.$}
\end{array} \right.$$
Then we set
\begin{equation} \label{defn}
Nf(x) = \sup_{0<r<R} \frac{1}{1+\|\vphi_{x,r,R}\|_{L^1(\mu)}}
\int |\vphi_{x,r,R}\, f|\, d\mu,
\end{equation}
for $f\in L^1_{loc}(\mu)$ and $x\in\supp(\mu)$.

Throughout all the paper $w$ stands for a positive weight in
$L^1_{loc}(\mu)$.
Sometimes the measure $w\,d\mu$ is denoted simply by $w$. The notation for
the dual weight is $\sigma := w^{-1/(p-1)}$, with $1<p<\infty$.

The first result that we will prove deals with the weak $(p,p)$ estimates.

\begin{theorem} \label{mainw1}
Let $p,\,\gamma$ be constants with $1\leq p<\infty$ and $0<\gamma\leq1$.
Let  $w$ be a positive weight. The following
statements are equivalent:
\begin{enumerate}
\item[(a)] All operators $T\in CZO(\gamma)$ are of weak
type $(p,p)$ with respect to $w\,d\mu$.

\item[(b)] The maximal operator $N$ is  of weak type $(p,p)$ with
respect to $w\,d\mu$.
\end{enumerate}
\end{theorem}

Next we state
the corresponding result for the strong $(p,p)$ estimates.

\begin{theorem} \label{main1}
Let $p,\,\gamma$ be constants with $1<p<\infty$ and $0<\gamma\leq1$.
Let  $w$ be a positive weight. The following
statements are equivalent:
\begin{enumerate}

\item[(a)] All operators $T\in CZO(\gamma)$ are bounded on $L^p(w)$.

\item[(b)] The maximal operator $N$ is bounded on $L^p(w)$.
\end{enumerate}
\end{theorem}

Let us denote by $Z_p$ the class of weights $w$ such that $N$
is bounded on $L^p(w)$, and by $Z_p^{weak}$ its weak version, that
is, the class of weights $w$ such that $N$ is
bounded from $L^p(w)$ into $L^{p,\infty}(w)$. Notice that
since $N$ is
bounded on $L^\infty(w)$, by interpolation we have $Z_p \subset Z_q$
if $1<p\leq q<\infty$.
On the other hand, the inclusion $Z_p\subset Z_p^{weak}$ is trivial, and by
duality (of CZO's) and Theorem \ref{main1} it follows that $w\in Z_p$ if and
only if $\sigma\in Z_{p'}$, where $p'$ stands for the conjugate exponent
of $p$, i.e. $p'=p/(p-1)$.

We will prove the following self improving property for this type weights:

\begin{theorem} \label{selfimprov}
Let $w$ be a positive weight and $1<p<\infty$.
If $w\in Z_p^{weak}$ and $\sigma=w^{-1/(p-1)} \in Z_{p'}^{weak}$,
then $w\in Z_p$ and $\sigma\in Z_{p'}$.
\end{theorem}

More detailed results are stated in Lemmas \ref{mainw} and
\ref{main} in Section 4.
In particular, necessary and sufficient conditions of ``Sawyer type''
are given for the boundedness of $N$ on $L^p(w)$ and also for the weak
$(p,p)$ case. Moreover, it is shown that if $w\in Z_p$ ($w\in
Z_p^{weak}$), then the maximal CZO
$$T_* f(x)= \sup_{\ve>0} |T_\ve f(x)|$$
is bounded on $L^p(w)$ [of weak type $(p,p)$ with respect to $w$].

Let us see an easy consequence of our results.
Given $\lambda\geq1$, let $M_{\lambda}$ be following version of the maximal
Hardy-Littlewood operator:
\begin{equation} \label{HLC}
M_\lambda f(x) = \sup_{r>0} \frac{1}{\mu(B(x,\lambda r))}
\int_{B(x,r)} |f|\,d\mu, \qquad \mbox{$f\in L^1_{loc}(\mu)$,
$x\in \supp(\mu)$}.
\end{equation}
It is easily seen that for any $\lambda\geq 1$,
\begin{equation} \label{mmm**}
Nf(x) \leq C(\lambda)\,M_\lambda f(x),\qquad
\mbox{$f\in L^1_{loc}(\mu)$,
$x\in \supp(\mu)$}.
\end{equation}
Thus all weights $w$ such that $M_\lambda$ is bounded on $L^p(w)$
belong to $Z_p$, and then all CZO's are bounded on  $L^p(w)$. In
particular, $A_p\subset Z_p$ if $1<p<\infty$.

Observe that the maximal operator $N$ is a centered maximal operator,
which is not equivalent to any ``reasonable'' non centered maximal
operator, as far as we know. This fact and the absence of any doubling
condition on $\mu$ are responsible for most of the difficulties that
arise in our arguments. For instance, it turns out that the weights of
the class $Z_p$ don't satisfy a reverse H\"older inequality, in
general. Indeed there are examples which show that it may happen that
$w\in Z_p$ but $w^{1+\ve}\not\in L^1_{loc}(\mu)$ for any $\ve>0$ (see
Examples \ref{ex3} and \ref{ex4}). Also, we will show that the weights
in $Z_p$
satisfy a  property much weaker than the $A_\infty$ condition of the
classical $A_p$ weights (see Definition \ref{sinfinit} and Lemma
\ref{cubspetits}), which is more difficult to deal with.

Let us notice that it has been shown in \cite{OP}
that, even with $\mu$ non doubling, if $w\in A_p$, then $w$ satisfies
a reverse H\"older inequality. As a consequence, $A_p\neq Z_p$ in general.

The plan of the paper is the following.
In Section 2 we show some examples which illustrate our results.
In Section 3 we recall the basic properties of the lattice of cubes
introduced in \cite{Tolsa3} and \cite{Tolsa4}, together with its
associated approximation of the identity. This construction will be
an essential tool for our arguments. In the same section we
will study some of the properties of the maximal operator $N$.
In Section 4 we state
Lemmas \ref{mainw} and \ref{main}, from which Theorems \ref{mainw1},
\ref{main1} and \ref{selfimprov} follow directly. Lemma \ref{mainw}
deals with the weak $(p,p)$ estimates, and it is proved in Sections
5--7, while the strong $(p,p)$ case is treated in Lemma \ref{main} and
is proved in Sections 8--10.
In Section 11 we explain how to prove the theorems above in their full
generality, without a technical assumption that is used in Sections 5--10
for simplicity. Finally, in Section 12 we show
which is the relationship between $Z_p$ and
$\rbmo(\mu)$ (this is the space of type $\bmo$ introduced in
\cite{Tolsa2}), and we make some remarks. In particular in this section we prove the following result:

\begin{theorem} \label{rbmo*}
Let $1<p<\infty$.
If $f\in\rbmo(\mu)$ and $\ve= \ve(\|f\|_*,p)>0$ is small enough,
then $e^{\ve f}\in Z_p$.
\end{theorem}

For the precise definition of $\rbmo(\mu)$, see Section 12.

% ***************************************************************************
% ***************************************************************************
% ***************************************************************************

\section{Some examples}

\begin{example} \label{ex1} If
$\mu(B(x,r))\approx r^n$ for all $x\in\supp(\mu)$,
  then $N f(x)\approx Mf(x)$, where $M$ is the usual centered
Hardy-Littlewood operator (defined in \rf{HLC} with $\lambda=1$). In
this case, the class $Z_p$ coincides with the class $A_p$.
\end{example}

\begin{example} \label{ex2}
In $\R^2$, consider the square $Q_0 = [0,1]^2$ and the measure $d\mu=
  \chi_{Q_0}dm$, where $dm$ stands for the planar Lebesgue
  measure, and take $n=1$. That is, we are interested in studying the
  weights for $1$-dimensional CZO's such as the Cauchy transform.
  Notice that $\mu$ is a doubling measure which does not satisfy the
  assumption in Example \ref{ex1}. For this measure,
  we have the uniform
  estimate $\int \frac{1}{|y-x|}\,d\mu(y)\leq C$.
  Then, from Theorem
  \ref{main1}, we deduce that the class $Z_p$ coincides with the class
  of $L^p$ weights for the fractional integral
$$I_0 f(x) = \int \frac1{|y-x|}\,f(y)\,d\mu(y),$$
since $Nf(x)\approx I_0|f|(x)$.
  This is the result that should be expected because, with our choice of
  $\mu$, $I_0$ is a CZO, and for all other $T\in CZO(\gamma)$, we have
  $|Tf(x)| \leq C_1\,I_0|f|(x)$.
\end{example}

\begin{example} \label{ex3}
This is an example studied by Saksman in his paper about weights
for the Hilbert transform \cite{Saksman}. We are in $\R$ and $n=1$.
Let $\ell_k=1/k!$ and
consider the intervals $I_k=\bigl(\frac1k - \frac{\ell_k}4,\,\frac1k +
\frac{\ell_k}4\bigr)$ for $k\geq 1$. Let $\mu$ be the Lebesgue measure
restricted to the set $S:=\bigcup_{k=1}^\infty I_k$.

Let $w$ be a
weight such that $w\geq1$ and $w_{|I_k}$ is constant for each $k\geq
1$. In \cite{Saksman}, it is proved that, for any $p\in(1,\infty)$, the
Hilbert transform is bounded
on $L^p(w)$ if and only if $w\in L^1(\mu)$. Almost the same
calculations show that the operators $S_k$ 
(defined after Lemma \ref{amesmes} below)
are uniformly bounded on
$L^p(w)$ if and only if $w\in L^1(\mu)$.

So, if a weight $w_0$ is defined by ${w_0}_{|I_k}=(n-2)!$, then
$w_0\in Z_p$ for all $p\in(1,\infty)$, by Lemma \ref{main} below. However, it is easily seen
that $w^{1+\ve}\not\in L^1(\mu)$ for any $\ve>0$. Therefore, $w_0$
does not satisfy a reverse H\"older inequality.
\end{example}

\begin{example} \label{ex4}
In this example we will show that there are measures $\mu$ and weights
$w\in Z_p$ such that the (centered) maximal Hardy-Littlewood operator $M$
is not bounded on $L^p(w)$. Also we will see that it may happen $w\in Z_p$
but $w\not\in Z_{p-\ve}$ for any $\ve>0$.

We take $d=n=1$. Suppose that $I_1$ and $I_2$ are disjoint intervals on
$\R$. The measure $\mu$ is the Lebesgue measure restricted to $I_1\cup I_2$.
Suppose that $\mu(I_1)=\mu(I_2)=L$, and let $D=\dist(I_1,I_2)$, with
$D\geq2L$. For $f=\chi_{I_1}$, the inequality $\|Mf\|_{L^p(w)} \leq
C_3\|f\|_{L^p(w)}$ implies
\begin{equation} \label{eqex41}
w(I_2) \leq C_4\,w(I_1),
\end{equation}
with $C_4$ depending on $C_3$ but not on $D$ or $L$.
By symmetry, \rf{eqex41} also holds interchanging $I_1$ and $I_2$.

Also, if $w\in Z_p$, from $\|N f\|_{L^p(w)} \leq
C_5\|f\|_{L^p(w)}$ we get
$\frac{L}{D}\,w(I_2)^{1/p} \leq C w(I_1)^{1/p}.$
That is,
\begin{equation} \label{eqex42}
\left(\frac{L}{D}\right)^p w(I_2) \leq C_6 w(I_1).
\end{equation}
The constant $C_6$ depends only on $C_5$.
By symmetry, we deduce
\begin{equation}\label{eqex43}
C_6^{-1} \left(\frac{L}{D}\right)^p w(I_2) \leq w(I_1)
\leq C_6 \left(\frac{L}{D}\right)^{-p} w(I_2)
\end{equation}
If $w$ is constant on each interval $I_1$, $I_2$, then $N$ is bounded
on $L^p(w)$ and it is easily seen that $\|N\|_{L^p(w)\rightarrow
L^p(w)}\leq C(C_6)$.

Now we introduce a new measure $\mu$ on $\R$. For each integer $m\geq1$
we consider the intervals $I_1^m = 1000^m + [-m-1,-m]$ and
$I_2^m = 1000^m + [m,m+1]$, so that $D_m := \dist(I_1^m,I_2^m) = 2m$ and $L=1$.
The measure $\mu$ is the Lebesgue measure restricted to
$\bigcup_{m=1}^\infty (I_1^m \cup I_2^m)$. The weight $w$ is constant
on each interval $w_{I_j^m}$, with $w_{|I_1^m} \equiv 1$ and
$w_{|I_2^m} \equiv (D_m/L)^p = (2m)^p$.

The maximal operator $M$ is not bounded on $L^p(w)$ because, otherwise,
we should have $w(I_2^m) \leq C\,w(I_1^m)$ uniformly on $m$, as
in \rf{eqex41}. On the other hand, \rf{eqex43} is satisfied (with the
corresponding subindices and superindices $m$) uniformly on $m$.
Taking also into account that $I_1^m\cup I_2^m$ is very far from
$I_1^r\cup I_2^r$ if $m\neq r$, it is easily checked that $N$ is bounded
on $L^p(w)$. Moreover, $N$ is not bounded on $L^{p-\ve}(w)$ for any $\ve>0$
because the inequality
$$\left(\frac{L}{D}\right)^{p-\ve} w(I_2) \leq C\,w(I_1)$$
fails for $m$ big enough.
\end{example}

% ***************************************************************************
% ***************************************************************************
% ***************************************************************************

\section{Preliminaries}

\subsection{The lattice of cubes}

For definiteness, by a cube we mean a closed cube with sides parallel
to the coordinate axes. We will assume that the constant $C_0$ in
\rf{creixlin} has been chosen big enough so that for all cubes $Q\subset
\R^d$ we have $\mu(Q) \leq C_0\, \ell(Q)^n$, where $\ell(Q)$ stands for the
side length of $Q$.

Given $\alpha,\beta>1$, we say that the cube $Q\subset\R^d$ is
$(\alpha,\beta)$-doubling if $\mu(\alpha Q) \leq \beta \mu(Q)$.
If $\alpha$ and $\beta$ are not specified and we say that some cube is
doubling, we are assuming $\alpha=2$ and $\beta$ equal to some
constant big enough ($\beta>2^d$, for example) which may depend from
the context.

\begin{rem} \label{moltdob}
Due to the fact that $\mu$ satisfies the growth condition \rf{creixlin}, there
are a lot ``big'' doubling cubes. To be precise, given any point
$x\in\supp(\mu)$ and $c>0$, there exists some
$(\alpha,\beta)$-doubling cube $Q$ centered at $x$ with $l(Q)\geq c$. This
follows easily from \rf{creixlin} and the fact that we are assuming
that $\beta>\alpha^n$.

On the other hand, if $\beta>\alpha^d$, then for $\mu$-a.e.\ $x\in\R^d$ there
exists a sequence of $(\alpha,\beta)$-doubling cubes $\{Q_k\}_k$ centered
at $x$ with $\ell(Q_k)\to0$ as $k\to \infty$.
So there are a lot of ``small'' doubling cubes too.
\end{rem}

Given cubes $Q,R$, with $Q\subset R$, we denote by $z_Q$ the center of $Q$, and by
$Q_R$ the smallest cube concentric with $Q$ containing $Q$ and
$R$. We set
$$\delta(Q,R) = \int_{Q_R\setminus Q}
\frac{1}{|x-z_Q|^n}\, d\mu(x).$$

We may treat points $x\in\R^d$ and the whole space $\R^d$ as if they were cubes (with
$\ell(x)=0$, $\ell(\R^d)=\infty$).
So for $x\in\R^d$ and some cube $Q$, the notations
$\delta(x,Q)$, $\delta(Q,\R^d)$ make sense. Of course, it may happen
$\delta(x,Q) = \infty$ and $\delta(Q,\R^d)=\infty$.

In the following lemma, proved in \cite{Tolsa3}, we recall some useful
properties of $\delta(\cdot,\cdot)$.

\begin{lemma}  \label{propdelta}

Let $P, Q, R\subset \R^d$ be cubes with  $P\subset Q\subset R$
The following properties hold:
\begin{itemize}
\item[(a)] If $\ell(Q) \approx \ell(R)$, then
$\delta(Q,R) \leq C$. In particular, $\delta(Q,\rho Q)\leq C_0\,2^n\,\rho^n$ for
$\rho>1$.

\item[(b)] If $Q\subset R$ are concentric and there are no doubling
cubes of the form $2^k Q$, $k\geq0$, with $Q\subset 2^k Q\subset R$, then
$\delta(Q,R)\leq C_7$.

\item[(c)]
$\ds \delta(Q,R) \leq C\,\left(1+\log\frac{\ell(R)}{\ell(Q)}\right).$

\item[(d)]
$\bigl|\delta(P,R) - [ \delta(P,Q) + \delta(Q,R) ] \bigr| \leq \ve_0.$
That is, with a different notation,
$\delta(P,R) = \delta(P,Q) + \delta(Q,R) \pm \ve_0$.

\end{itemize}
\end{lemma}

The constants $C$ and $\ve_0$ that appear in (b), (c) and (d) depend
on $C_0,n,d$.
The constant $C$ in (a) depends, further, on the constants that are
implicit in the relation $\approx$.
Let us insist on the fact that a notation such as $a= b \pm \ve$ does not mean
any precise equality but the estimate $|a-b|\leq \ve$.

Now we will describe the lattice of cubes introduced in
\cite{Tolsa4}. In the following lemma, $Q_{x,k}$ stands for a cube
centered at $x$, and we allow $Q_{x,k}=x$ and  $Q_{x,k}=\R^d$. If
$Q_{x,k}\neq x,\R^d$, we say that $Q_{x,k}$ is a {\em transit
cube}.

\begin{lemma} \label{lattice}
Let $A$ be an arbitrary positive constant big enough. There exists a
family of cubes $Q_{x,k}$, for all $x\in\supp(\mu)$, $k\in\Z$,
centered at $x$, and such that:
\begin{itemize}
\item[(a)] $Q_{x,k} \subset Q_{x,j}$ if $k\geq j$.

\item[(b)] $\lim_{k\to+\infty} \ell(Q_{x,k})= 0$ and
 $\lim_{k\to-\infty} \ell(Q_{x,k})= \infty$.

\item[(c)] $\delta(Q_{x,k},Q_{x,j}) = (j-k)A \pm \ve$ if $j>k$ and
  $Q_{x,k},\,Q_{x,j}$ are transit cubes.

\item[(d)] $\delta(Q_{x,k},Q_{x,j}) \leq (j-k)A + \ve$ if $j>k$.

\item[(e)] If $2Q_{x,k} \cap 2Q_{y,k}\neq \varnothing$, then $2Q_{x,k}
 \subset Q_{y,k-1}$ and $\ell(Q_{x,k}) \leq \ell(Q_{y,k-1})/100$.

\item[(f)] There exists some $\eta>0$ such that if $m\geq1$ and
  $2Q_{x,k+m} \cap 2Q_{y,k}\neq \varnothing$, then
  $\ell(Q_{x,k+m}) \leq 2^{-\eta A m} \ell(Q_{y,k})$.

%\item[(g)] $Q_{x,k}$ is doubling.
\end{itemize}

The constants $\ve,\,\eta$ in (c), (d) and (f) depend
on $C_0,n,d$, but not on $A$.
\end{lemma}

See \cite[Section 3]{Tolsa4} for the proof.
The constant $\ve$ above must be understood as an error
term, because we will take $A\gg \ve$.
Let us notice also that, if
necessary, the cubes $Q_{x,k}$ can be chosen so that they are doubling
(see \cite{Tolsa4}). However we don't need this assumption.

\begin{rem}
If $x\in\supp(\mu)$ is such that $\int_{B(x,1)} |y-x|^{-n}\,d\mu(y)
<\infty$, then it follows from the properties of the lattice that
there exists some $K_x\in\Z$ such that $Q_{x,k}=x$ for $k>K_x$ and
$Q_{x,k}\neq x$ for $k\leq K_x$. In this case we say that $Q_{x,k}$ is a
{\em stopping cube} (or {\em stopping point}).

If $\int_{\R^d\setminus B(x,1)} |y-x|^{-n}\,d\mu(y)<\infty$ (which
does not depend on $x\in\supp(\mu)$), then there exists some constant
$\bar K_x$ such that $Q_{x,k} = \R^d$ for $k<\bar K_x$ and
$Q_{x,k} \neq \R^d$ for $k\geq\bar K_x$.
We say that $\R^d$ is an (or the) {\em initial cube}.
From the property (e) in
the lemma above, it follows easily that $|\bar K_x - \bar
K_y|\leq 1$ for $x,y\in\supp(\mu)$. However, as shown in
\cite{Tolsa4}, the construction of the
lattice can be done so that $\bar K_x = \bar K_y=:\bar
K_0$ for all $x,y$, and so that $\delta(Q_{x,\bar K_0+m},\R^d) =
mA\pm\ve$ for $m\geq1$.
For simplicity, we will assume that
our lattice fulfils these properties.

If $\int_{B(x,1)} |y-x|^{-n}\,d\mu(y) = \int_{\R^d\setminus B(x,1)}
|y-x|^{-n}\,d\mu(y) = \infty$, then all the cubes $Q_{x,k}$, $k\in\Z$,
satisfy $0<\ell(Q_{x,k})<\infty$. That is, they are transit cubes.

\end{rem}

We denote $\DD_k=\{Q_{x,k}:\,x\in\supp(\mu)\}$ for $k\in\Z$, and
$\DD=\bigcup_{k\in\Z}\DD_k$.

Consider a cube $Q\subset\R^d$ whose center may not be in
$\supp(\mu)$. Let $Q_{x,k}$ be one of the smallest cubes in $\DD$
containing $Q$ in the following sense. Set
$$\ell=\inf\{\ell(Q_{x,j}):\,
Q_{x,j}\in\DD,\, Q\subset Q_{x,j}\}.$$
Take $Q_{x,k}$ containing $Q$ such that $\ell(Q_{x,k})\leq
\frac{100}{99} \ell$.  Then we write $Q\in \AD_k$ (by the property (e)
in Lemma \ref{lattice},
k depends only on $Q$). In a sense, $Q$ is approximately in $\DD_k$.
Given $k,j$ with $-\infty\leq k\leq j \leq +\infty$, we also denote
$\AD_{k,j} = \bigcup_{h=k}^j \AD_h$. If $Q$ is such that there are
cubes $Q_{x,k}$,  $Q_{y,k-1}$ with $Q_{x,k}\subset Q \subset
Q_{y,k-1}$, then it follows easily that $Q\in\AD_{k,k-1}$.

% ***************************************************************************
% ***************************************************************************
% ***************************************************************************

\subsection{The kernels $s_k(x,y)$}

For each $x\in\supp(\mu)$, $s_k(x,\cdot)$ is a non negative radial non
increasing function with center $x$, supported on $2Q_{x,k-1}$, and such that
\begin{itemize}
\item[(a)] $\ds s_k(x,y) \leq \frac{1}{A|x-y|^n}$ for all $y\in\R^d$.
\item[(b)] $\ds s_k(x,y) \approx  \frac{1}{A\ell(Q_{x,k})^n}$ for all $y\in
Q_{x,k}$.
\item [(c)] $\ds s_k(x,y) =  \frac{1}{A|x-y|^n}$ for all $y\in
  Q_{x,k-1}\setminus Q_{x,k}$.
\item[(d)] $\ds \nabla_{\!y} s_k(x,y) \leq C\,A^{-1}\,\min\left(
 \frac1{\ell(Q_{x,k})^{n+1}},\, \frac1{|x-y|^{n+1}} \right)$ for all $y\in\R^d$.
\end{itemize}

\begin{lemma} \label{suportpetit}
If $y\in\supp(\mu)$, then $\supp(s_k(\cdot,y)) \subset Q_{y,k-2}$.
If $Q\in\AD_k$ and $z\in Q\cap\supp(\mu)$, then $\supp(s_{k+m}(z,\cdot))\subset
\frac{11}{10}Q$ for all $m\geq3$, and  $\supp(s_{k+m}(\cdot,z))\subset
\frac{11}{10}Q$ for all $m\geq4$.
\end{lemma}

\begin{proof}
For the assertion $\supp(s_k(\cdot,y)) \subset Q_{y,k-2}$, see
\cite{Tolsa3} or \cite{Tolsa4}.

Let $Q\in\AD_k$ and $z\in Q\cap\supp \mu$.
We have $Q\not\subset Q_{z,k+1}$, because otherwise $Q\not\in
\AD_k$. Thus $\ell(Q_{z,k+1}) \leq 2\ell(Q)$. Then,
$$\supp(s_{k+m}(z,\cdot)) \subset 2Q_{z,k+m-1} \subset\frac{11}{10}Q,$$
because
$\ell(2Q_{z,k+m-1})\leq\frac{2}{100}\ell(Q_{z,k+1})\leq\frac{4}{100}\ell(Q)$.
Finally, the inclusion $\supp(s_{k+m}(\cdot,z))\subset
\frac{11}{10}Q$ follows in a similar way.
\end{proof}

In \cite[Section 3]{Tolsa4} the following estimates are proved.

\begin{lemma}  \label{hoho}
If $A$ is big enough, then for all $k\in\Z$ and $z\in\supp(\mu)$ we have
\begin{equation} \label{mmhh1}
\int s_k(z,y)\,d\mu(y) \leq \frac{10}{9}\qquad \mbox{and} \qquad
\int s_k(x,z)\,d\mu(x) \leq \frac{10}{9}.
\end{equation}
If moreover $Q_{z,k}$ is a transit cube, then
\begin{equation} \label{mmhh2}
\int s_k(z,y)\,d\mu(y) \geq \frac{9}{10}\qquad \mbox{and} \qquad
\int s_k(x,z)\,d\mu(x) \geq \frac{9}{10}.
\end{equation}
\end{lemma}

In the following lemma we state
another technical result that we will need.

\begin{lemma} \label{amesmes}
For all $k\in\Z$ and $x,y\in\supp(\mu)$, we have
\begin{equation} \label{quasiad}
s_k(x,y) \leq C \bigl(s_{k-1}(y,x) + s_k(y,x) + s_{k+1}(y,x)\bigr).
\end{equation}
\end{lemma}

The proof follows easily from our construction. See also
\cite[Lemma 7.8]{Tolsa3}.

We will denote by $S_k$ the integral operator associated with the
kernel $s_k(x,y)$ and the measure $\mu$.
Observe that \rf{mmhh1} implies that
the operators $S_k$ are bounded uniformly on
$L^p(\mu)$, for all $p\in[1,\infty]$. Also, from \rf{quasiad} we get
\begin{equation} \label{quho}
S_kf(x) \leq C\bigl(
S_{k-1}^*f(x) + S_{k}^*f(x) + S_{k+1}^*f(x)\bigr),
\end{equation}
for $f\in L^1_{loc}(\mu)$, $f\geq0$, and $x\in\R^d$.

Notice that the only property in the definition of Calder\'on-Zygmund
kernel which $s_k(x,y)$ may not satisfy is the gradient condition on the first
variable. It is not difficult to check that if
the functions $\ell(Q_{x,k})$
were Lipschitz (with respect to $x$) uniformly on $k$, then
we would be able to define  $s_k(x,y)$ so that
\begin{equation}\label{regx}
|\nabla_{\!x} s_k(x,y)(y)| \leq \frac{C}{A|x-y|^{n+1}},
\end{equation}
in addition to the properties above. The following lemma solves
this question.

\begin{lemma} \label{SCZO}
The lattice $\DD$ can be constructed so that
the functions $\ell(Q_{x,k})$
are Lipschitz (with respect to $x\in\supp(\mu)$) uniformly on $k$
and the properties (a)--(f) in Lemma
\ref{lattice} still hold. In this case, the operators $S_k$,
$k\in\Z$, are CZO's with constants uniform on $k$.
\end{lemma}

\begin{proof} Suppose that the cubes $Q_{x,k}\in\DD$ have already been
  chosen and the properties stated in Lemma \ref{lattice} hold.
Let us see how we can choose cubes $\wt{Q}_{x,k}$, substitutes of
$Q_{x,k}$, such that $\psi_k(x) := \ell(\wt{Q}_{x,k})$ are Lipschitz
functions on $\supp(\mu)$. For a fixed $k$, we set
\begin{equation} \label{noucubs}
\psi_k(x) := \sup_{z\in\supp(\mu)} (\ell(Q_{z,k}) - |x-z|).
\end{equation}
It is easily seen that this is a non negative Lipschitz function, with constant
independent of $k$. Then, we denote by $\wt{Q}_{x,k}$ the cube centered at $x$
with side length $\psi_k(x)$.

We have to show that $\wt{Q}_{x,k}$ is a good choice as a cube of the
scale $k$. Indeed, by \rf{noucubs} it is clear that
$\ell(\wt{Q}_{x,k}) \geq \ell(Q_{x,k})$. Thus
$Q_{x,k}\subset \wt{Q}_{x,k}.$

Take now $z_0\in\supp(\mu)$ such that
$$\ell(Q_{z_0,k}) - |x-z_0| \geq  \frac{99}{100}\,\ell(\wt{Q}_{x,k}).$$
We derive $|x-z_0|\leq \ell(Q_{z_0,k})$, and also
$\ell(\wt{Q}_{x,k})\leq 100\,\ell(Q_{z_0,k})/99$. Thus
$x\in 2Q_{z_0,k}$ and $\wt{Q}_{x,k} \subset 4Q_{z_0,k}$.
The inclusions
$Q_{x,k}\subset \wt{Q}_{x,k}\subset 4Q_{z_0,k}$
imply $\delta(\wt{Q}_{x,k}, Q_{x,k}) \leq C_8 \ll
\delta(\wt{Q}_{x,k}, Q_{x,k-1})$, with $C_8$ depending only on $n$,
$d$, $C_0$.
One can verfy that the properties in Lemma \ref{lattice} still hold,
assuming that the constant $C_8$ is absorbed by the ``error''
$\ve$ in (c) and (d) in Lemma \ref{lattice}.
\end{proof}

% ***************************************************************************
% ***************************************************************************
% ***************************************************************************

\subsection{The maximal operator $\boldsymbol N$}

In the following lemma we show which is the relationship between $N$
and the operators $S_k$.

\begin{lemma}
For all $f\in L^1_{loc}(\mu)$, $x\in \R^d$, we have
$$Nf(x) \approx \sup_{k\in\Z} S_k|f|(x),$$
with constants depending on $A,\,C_0,\,n,\,d$ but independent of
$f$ and $x$.
\end{lemma}

\begin{proof}
For fixed $x\in\supp(\mu)$ and $k\in\Z$, we have
$s_k(x,y) \leq C\,\vphi_{x,r,R}(y)$,
with $r=C\ell(Q_{x,k})$ and $R=C\ell(Q_{x,k-1})$. Assume $0<r,R<\infty$.
Since $\|\vphi_{x,r,R}\|_{L^1(\mu)} \leq C$ we get
$$S_k|f|(x) \leq \frac{C}{1+\|\vphi_{x,r,R}\|_{L^1(\mu)}}
\int|\vphi_{x,r,R}f|\,d\mu \leq C\, Nf(x).$$
If $r=0$ or $R=\infty$, we also have $S_k|f|(x) \leq C\, Nf(x)$
by an approximation argument, and so
$\sup_{k} S_k|f|(x) \leq C \,Nf(x).$

Let us see the converse inequality. Given $0<r<R<\infty$,
let $k$ be the least integer
such that $Q_{x,k}\subset B(x,r)$. Now let $m$ be the least positive integer
such that $B(x,R) \subset Q_{x,k-m}$.
Then we have
$$\vphi_{x,r,R}(y) \leq C\,(s_k(x,y) + s_{k-1}(x,y) + \cdots +
s_{k-m}(x,y)).$$
Also, it is easily checked that
$1+\|\vphi_{x,r,R}\|_{L^1(\mu)} \geq C^{-1}\,m.$
Therefore,
$$\frac{1}{1+\|\vphi_{x,r,R}\|_{L^1(\mu)}}
\int|\vphi_{x,r,R}f|\,d\mu
 \leq \frac{C}{m} \sum_{h=0}^m S_h|f|(x)
\leq  C \sup_i S_i|f|(x).$$
\end{proof}

In the rest of the paper we will assume that {\bf $N$ is
defined not by \rf{defn}, but as}
$$Nf(x) := \sup_{k\in\Z} S_k|f|(x).$$
With this new definition we have:

\begin{lemma} \label{open}
Let $\lambda>0$ and $f\in L^1_{loc}(\mu)$. For each $k\in\Z$, the set
$\{x\in\R^d:\,S_k|f|(x)>\lambda\}$ is open. As a consequence,
$\{x\in\R^d:\,Nf(x)>\lambda\}$ is open too.
\end{lemma}

The proof is an easy exercise which is left for the reader.

Given a fixed $x\in\supp(\mu)$, we can think of $S_k f(x)$ as an
average of the means
$m_{B(x,r)}f:=\int_{B(x,r)}f\,d\mu/\mu(B(x,r))$ over some range of
radii $r$. Arguing in this way,
\rf{mmm**} follows. We will exploit the same idea in the
following lemma.

\begin{lemma} \label{guai}
For all $\alpha>1$, we can choose constants $A,\,\beta,\,C_9$ big enough so
that the following property holds: Let $x\in\supp(\mu)$, $k\in\Z$ and
$f\in L^1_{loc}(\mu)$, and assume that $Q_{x,k}$ is a transit
cube. Then there exists some ball $B(x,r)$ with
$Q_{x,k}\subset B(x,\alpha^{-1}r),\,B(x,r) \subset Q_{x,k-1}$ such that
$B(x,\alpha^{-1}r)$ is $(\alpha,\beta)$-doubling and $m_{B(x,r)}|f| \leq
C_9\,S_k |f|(x)$.
\end{lemma}

It is easy to check that there are balls $B(x,r_1)$ and $B(x,r_2)$ with
$Q_{x,k}\subset B(x,\alpha^{-1}\,r_1),\, B(x,r_2)\subset Q_{x,k}$ such that
$B(x,\alpha^{-1}\,r_1)$ is $(\alpha,\beta)$-doubling
and $m_{B(x,r_2)} |f|\leq C\,S_k |f|(x)$. However, it is
more difficult to see that we may take $B(x,r_1)=B(x,r_2)$, as the lemma
asserts.

On the other hand, the lemma is false if we substitute the condition
``$m_{B(x,r)}|f| \leq C_9\,S_k |f|(x)$'' by
``$m_{B(x,r)}|f| \geq C_9^{-1}\,S_k |f|(x)$''.

\begin{proof} We denote $\lambda := S_k |f|(x)$, $R_0 = d^{1/2}
\ell(Q_{x,k})$,
and $R_1 = \ell(Q_{x,k-1})/2$.
Recall that, for fixed $x$, $k$, we have defined
$s_k(x,y)=\psi(|y-x|)$, where $\psi:\R\longrightarrow\R$ is non negative,
smooth, radial, and non increasing. Then,
$$\lambda = \int |f(y)| s_k(x,y) d\mu(y)
= \int_0^\infty |\psi'(r)| \left(\int_{B(x,r)}  |f|\,d\mu\right) dr.$$
We denote $h(r) = |\psi'(r)|\,\mu(B(x,\alpha^{-1}r))$ and
$$m_\alpha(r) = \frac{1}{\mu(B(x,\alpha^{-1}r))} \int_{B(x,r)}
|f|\,d\mu.$$
Thus, $\lambda = \int_0^\infty h(r)\, m_\alpha(r)\,dr.$

Let us see that $\int_{\alpha R_0}^{R_1} h(r)\,dr$ is big.
Recall that $\psi'(r)=1/(A\,r^{n+1})$ for $r\in [R_0, R_1]$.
Then, for $s\in [R_0,\,\alpha^{-1}R_1]$, we have
$|\psi'(\alpha s)| = |\psi'(s)|/\alpha^{n+1}$. Therefore,
\begin{eqnarray*}
\int_{\alpha R_0}^{R_1} h(r)\,dr & = &\alpha^{-n} \int_{R_0}^{\alpha^{-1}R_1}
|\psi'(s)|\,\mu(B(x,s))\,ds  \\
& = &  \alpha^{-n} \Biggl(
\int_{R_0\leq |x-y| \leq \alpha^{-1}R_1} s_k(x,y)\,d\mu(y) \\
& & \mbox{} + \psi(R_0)\,\mu(B(x,R_0)) -
\psi(\alpha^{-1} R_1)\,\mu(B(x,\alpha^{-1}R_1)) \Biggr)\\
& \geq & \alpha^{-n} \left(\int_{R_0\leq |x-y| \leq \alpha^{-1}R_1}
s_k(x,y)\,d\mu(y) - C_0\,A^{-1}\right).
\end{eqnarray*}
Since $\int_{|x-y|\leq R_0} s_k(x,y)\,d\mu(y) \leq C\,A^{-1}$ and
also  $\int_{|x-y|\geq \alpha^{-1}R_1} s_k(x,y)\,d\mu(y) \leq C\,A^{-1}$,
for $A$ big enough we obtain
$$\int_{\alpha R_0}^{R_1} h(r)\,dr \geq \frac{1}{2\alpha^n}=:M,$$
using \rf{mmhh2}.
If we denote the measure $h(r)\,dr$ by $h$, we get
$$h\{r\geq0:m_\alpha(r)>2\lambda/M\}
\leq \frac{M}{2\lambda} \int_0^\infty m_\alpha(r)\, h(r)
\,dr = \frac{M}{2}.$$
Thus,
$$h\{r\in[\alpha R_0,R_1]:\, m_\alpha(r) \leq 2\lambda/M\} \geq M - \frac{M}{2} =
\frac{M}{2}.$$

Now we will deal with the doubling property. If $B(x,\alpha^{-1}r)$ is
not $(\alpha,\beta)$-doubling, we write $r\in ND$. We have
\begin{eqnarray*}
h([\alpha R_0,\,R_1]\cap ND) & = & \int_{r\in [\alpha R_0,R_1]\cap ND}
|\psi'(r)|\, \mu(B(x,\alpha^{-1}r)) \,dr \\
& \leq & \beta^{-1} \int_{\alpha R_0}^{R_1} |\psi'(r)|\,\mu(B(x,r)) \,dr \\
& \leq & \beta^{-1}\, \int s_k(x,y)\,d\mu(y) \leq \frac{10}9\beta^{-1}.
\end{eqnarray*}
Therefore,
$$h\bigl( \{r\in [\alpha R_0,\,R_1]:\, m_\alpha(r) \leq
2 \lambda/M\} \setminus
ND\bigr) \geq \frac{M}{2} - \frac{10}9\beta^{-1}.$$

So if we take $\beta$ big enough, there exists some $r \in[\alpha R_0,R_1]$
such that $B(x,\alpha^{-1}r)$ is $(\alpha,
\beta)$-doubling and $m_{B(x,r)}|f| \leq
m_\alpha(r) \leq 2\lambda/M$.
\end{proof}

As a direct corollary of Lemma \ref{guai} we get:

\begin{lemma} \label{guai'}
Assume that $A,\,\beta,\,C_{10}$ are positive and big enough. Let
$x\in\supp(\mu)$, $k\in\Z$ and
$f\in L^1_{loc}(\mu)$. If $Q_{x,k}$ is a transit cube, then there
exists some $(2,\beta)$-doubling cube
$Q\in \AD_{k,k-1}$ centered at $x$ such that
$m_{2Q}|f| \leq C_{10}\,S_k |f|(x)$.
\end{lemma}

In the rest of the paper we will assume that the constant $A$ used to
construct the lattice $\DD$ and the kernels $s_k(x,y)$ has been chosen big
enough so that the conclusion of the preceding lemma holds.

% ***************************************************************************
% ***************************************************************************
% ***************************************************************************

\section{The main lemmas}

Theorems \ref{mainw1}, \ref{main1} and \ref{selfimprov} follow from the
following two lemmas:

\begin{lemma} \label{mainw}
Let $p,\,\gamma$ be constants with $1\leq p<\infty$ and $0<\gamma\leq1$.
Let  $w>0$ be a weight and $\sigma = w^{-1/(p-1)}$ (for $p\neq1$). 
The following statements are equivalent:
\begin{enumerate}

\item[(a)] All operators $T\in CZO(\gamma)$ are of weak
type $(p,p)$ with respect to $w\,d\mu$.

\item[(b)] For all $T\in CZO(\gamma)$, $T_*$ is of weak type $(p,p)$ with
respect to $w\,d\mu$.

\item[(c)] The maximal operator $N$ is  of weak type $(p,p)$ with
respect to $w\,d\mu$.

\item[(d)] The operators $S_k$ are of weak type $(p,p)$ with respect
to $w\,d\mu$ uniformly on $k\in\Z$.

\item[(e)] (Only in the case $p\neq1$). For all $k\in\Z$ and all cubes $Q$,
\begin{equation} \label{cond1'}
\int |S_k(w\,\chi_Q)|^{p'}\, \sigma\,d\mu \leq C\,w(Q),
\end{equation}
with $C$ independent of $k$ and $Q$.
\end{enumerate}
\end{lemma}

\begin{lemma} \label{main}
Let $p,\,\gamma$ be constants with $1<p<\infty$ and $0<\gamma\leq1$.
Let  $w>0$ be a weight and $\sigma = w^{-1/(p-1)}$. The following
statements are equivalent:
\begin{enumerate}

\item[(a)] All operators $T\in CZO(\gamma)$ are bounded on $L^p(w)$.

\item[(b)] For all $T\in CZO(\gamma)$, $T_*$ is bounded on $L^p(w)$.

\item[(c)] The maximal operator $N$ is bounded on $L^p(w)$.

\item[(d)] The operators $S_k$ are bounded on $L^p(w)$ uniformly on $k\in\Z$.

\item[(e)] For all $k\in\Z$ and all cubes $Q$,
$$\int |S_k(\sigma\,\chi_Q)|^p\, w\,d\mu \leq C\,\sigma(Q)$$
and
$$\int |S_k(w\,\chi_Q)|^{p'}\, \sigma\,d\mu \leq C\,w(Q),$$
with $C$ independent of $k$ and $Q$.
\end{enumerate}

\end{lemma}

%Let us remark that if we assume $w\geq0$ instead of $w>0$, any of the
%statements (a)-(d) implies that either $w>0$, or $w\equiv0$, or
%$w$ is identically infinite.

Notice that the Sawyer type conditions (e) in Lemma \ref{mainw} and
Lemma \ref{main} involve the operators $S_k$ instead of the maximal
operator $N$. In the present formulation these conditions are much
weaker and of more geometric nature than the analogous conditions
involving $N$.

The scheme for proving both lemmas is the same.
In both cases we will start by (c) $\Rightarrow$ (b).
Later we will see  (b) $\Rightarrow$ (a) $\Rightarrow$ (d) $\Rightarrow$ (e).
These will be the easy implications. Notice, for instance, that (b)
$\Rightarrow$ (a) is trivial.
Finally we will show (e) $\Rightarrow$ (c) (except in the weak $(1,1)$
case). This will be the most
difficult part of the proof (in both lemmas). In the weak
$(1,1)$ case, we will see directly the implication (d) $\Rightarrow$ (c).

For simplicity, to prove Lemmas \ref{mainw} and \ref{main}, we will
assume that {\bf all the cubes $Q_{x,k}\in\DD$ are transit cubes}. In
Section 11 we will give some hints for the proof in the general case.
We have operated in this way because the presence of stopping cubes in
the lattice $\DD$
introduces some technical difficulties which make the proofs more
lengthy, but the ideas and arguments involved are basically the same
than in the special case in which all the cubes in $\DD$ are transit cubes.

First we will prove Lemma \ref{mainw}.

% ***************************************************************************
% ***************************************************************************
% ***************************************************************************

\section{The implication {\rm (c)} $\Rightarrow$ {\rm (b)} of Lemma
\ref{mainw}}

\begin{definition}  \label{sinfinit}
We say that $w$ satisfies the $Z_\infty$ property if there exists some
constant
$\tau>0$ such that for any cube $Q\in \AD_{k}$ and any set
$A\subset \R^d$, if
\begin{equation}  \label{hsinfinit}
S_{k+3}\chi_A(x)\geq 1/4\qquad \mbox{for all $x\in Q$,}
\end{equation}
then $w(A\cap 2Q) \geq \tau\, w(Q)$.
\end{definition}

\begin{lemma} \label{lsinfinit}
If
$$\int |S_k^*(w\,\chi_Q)|^{p'}\, \sigma\,d\mu \leq C\,w(Q)$$
for all cubes $Q$ and all $k\in\Z$, then $w$ satisfies the property
$Z_\infty$.
\end{lemma}

\begin{proof}
Take $Q\in \AD_{k}$ and
$A\subset\R^d$ satisfying \rf{hsinfinit}.
By the assumption above, the fact that $\supp(s_{k+3}(x,\cdot))\subset
2Q$ for $x\in Q$, and H\"older's inequality, we get
\begin{eqnarray*}
w(Q) & \leq & 4 \int_Q (S_{k+3}\chi_A)\, w\,d\mu
 \\ & = & 4 \int (S_{k+3}\chi_{A\cap2Q})\, w\,d\mu \,
= \, 4 \int_{A\cap2Q} S_{k+3}^*(w\,\chi_Q)\,d\mu \\
& \leq & 4
\left(\int S_{k+3}^*(w\,\chi_Q)^{p'}\, \sigma\,d\mu\right)^{1/p'}\,
w(A\cap 2Q)^{1/p} \\
& \leq & C\, w(Q)^{1/p'}\, w(A\cap 2Q)^{1/p},
\end{eqnarray*}
and so $w(Q)\leq C\,w(A\cap 2Q)$.
\end{proof}

Notice that if $N$ is bounded on $L^p(w)$ for some $p\in(1,\infty)$, 
then the operators
$S_k$ are bounded on $L^p(w)$ uniformly on $k$. By duality, the operators
$S_k^*$ are bounded on $L^{p'}(\sigma)$ uniformly on $k$ too. Then, by Lemma
\ref{lsinfinit}, $w$ satisfies $Z_\infty$.

Occasionally we will apply the $Z_\infty$ condition by means of the
following lemma.

\begin{lemma}  \label{cubspetits}
Suppose that $w$ satisfies the $Z_\infty$ property.
Let $A\subset\R^d$ and $Q\in\AD_{h}$. Let $\{P_i\}_i$ be
a family of cubes with finite overlap such that $A\cap\frac{3}{2}Q\subset
\bigcup_i P_i$, with $P_i\in \AD_{+\infty,h+4}$ for all $i$. There
exists some constant $\delta>0$ such that if
$\mu(A \cap P_i) \leq \delta \,\mu(P_i)$ for each $i$, then
\begin{equation} \label{bm0.5}
w(2Q\setminus A) \geq \tau\,w(Q),
\end{equation}
for some constant $\tau>0$ (depending on $Z_\infty$).
If, moreover, $w(2Q) \leq C_{11}\,w(Q)$, then
\begin{equation} \label{bm1}
w(A\cap 2Q) \leq (1-C_{11}^{-1}\tau)\,w(2Q).
\end{equation}
\end{lemma}

This lemma, specially the inequality \rf{bm1}, shows that
the $Z_\infty$ property
can be considered as a weak version of the usual $A_\infty$ property
satisfied by the $A_p$ weights. Notice that unlike $A_\infty$, the
$Z_\infty$ condition is not symmetric on $\mu$ and $w$.

Let us remark that
we have not been able to prove that the constant $1-C_{11}^{-1}\tau$
in \rf{bm1} can be substituted by some constant $C_\delta$ tending to 0 as
$\delta\to0$. Many difficulties in the arguments below stem
from this fact.

\begin{proof}
For a fixed $x_0\in Q$,
we denote $Q_0 := Q_{x_0,h+3}$. Observe that
$$\supp(s_{h+3}(x_0,\cdot)\,\chi_A)\subset A\cap 2Q_{x_0,h+2}\subset
A\cap\frac{3}{2}Q \subset \bigcup_i P_i.$$
We have
$$S_{h+3} \chi_A(x_0) \leq C \sum_{j=1}^{n_0} \frac{\mu(2^jQ_0\cap
A)}{\ell(2^jQ_0)^n},$$
where $n_0$ is the least integer such that $2Q_{x_0,h+2}\subset
2^{n_0}Q_0$.
If $P_i\cap 2^jQ_0 \neq \varnothing$, then $\ell(P_i)\leq
\ell(2^jQ_0)/10$. Therefore,
$$\mu(2^jQ_0\cap A) \leq \!\!\sum_{i:\,P_i\cap 2^jQ_0\neq\varnothing}\!\!
\mu(P_i\cap A) \leq \!\! \sum_{i:\,P_i\cap 2^jQ_0\neq\varnothing}\!\! \delta\,
\mu(P_i) \leq C\,\delta \mu(2^{j+1}Q_0).$$
Therefore,
$$S_{h+3} \chi_A(x_0) \leq C\,\delta \sum_{j=1}^{n_0}
\frac{\mu(2^{j+1}Q_0)}{\ell(2^jQ_0)^n} \leq C\,\delta.$$
If $\delta$ is small enough, we have
$S_{h+3} \chi_{\R^d\setminus A}(x_0)\geq 1/4$ for all $x_0\in Q$,
by\rf{mmhh2}. Thus \rf{bm0.5} holds.

Finally, \rf{bm1} follows easily from \rf{bm0.5}.
\end{proof}

\vspace{5mm}
The implication {\rm (c)} $\Rightarrow$ {\rm (b)} of Lemma
\ref{mainw} is a direct consequence of
the good $\lambda$ inequality in next lemma.

\begin{lemma}  \label{lemagli}
Let $T\in CZO(\gamma)$ and $w$ which satisfies the $Z_\infty$ condition.
There exists some $\eta>0$ such that for all $\lambda,\ve>0$
\begin{equation} \label{gli1}
w\{x:\,T_*f(x)>(1+\ve)\,\lambda,\, Nf(x)\leq\delta\lambda\} \leq (1-\eta)
\,w\{x:\,T_*f(x) > \lambda\}
\end{equation}
if $\delta=\delta(\eta,\ve)>0$ is small enough.
\end{lemma}

The constant $\delta$ depends also on the weak (1,1) norm
of $T_*$ (with respect to $\mu$) and on $n,d$, besides of $\eta,\ve$,
but not on $\lambda$.

\begin{proof}
Given $\lambda>0$, we set $\Omega_\lambda = \{x:\,T_*f(x) > \lambda\}$ and
$$A_\lambda = \{x:\,T_*f(x)>(1+\ve)\,\lambda,\, Nf(x)\leq\delta\lambda\}.$$
So we have to see that there exists some $\eta>0$ such that, for all
$\ve>0$ and $\lambda>0$,
$w(A_\lambda)  \leq (1-\eta)\,w(\Omega_\lambda)$
if we choose $\delta=\delta(\eta,\ve)>0$ small enough.

Since $\Omega_\lambda$ is open, we can consider a Whitney decomposition of
it. That is, we set $\Omega_\lambda = \bigcup_i Q_i$, so that the cubes $Q_i$
have disjoint interiors, $\dist(Q_i,\R^d\setminus\Omega_\lambda)\approx
\ell(Q_i)$ for each $i$, and the cubes
$4Q_i$ have finite overlap.

Take a cube $Q_i$ such that there exists some $x_0\in 2Q_i$ with $Nf(x_0)\leq
\delta\lambda$. By standard arguments, one can check that for any $x\in 2Q_i$,
$$T_* (f\chi_{\R^d\setminus3Q_i})(x) \leq \lambda + C\, M_{R}^c f(x),$$
where $M_{R}^c$ is the centered radial maximal Hardy-Littlewood operator:
$$M_{R}^c f(x) = \sup_{r>0} \frac{1}{r^n} \int_{B(x,r)} |f|\, d\mu.$$
Since $M_{R}^c f\leq
C\,Nf$, we get $T_* f\chi_{\R^d\setminus3Q_i}(x) \leq (1+C\,\delta)\,\lambda$
if $x\in A_\lambda \cap 2Q_i$. For $\delta$ small enough, this implies
$T_* (f\chi_{3Q_i})(x) \geq \frac{\ve}{2}\, \lambda$
for all $x\in A_\lambda \cap 2Q_i$.
So we have
\begin{equation} \label{nn*1}
A_\lambda \cap 2Q_i\subset
\{x\in 2Q_i:\, T_*(f\chi_{3Q_i})(x)>\ve\lambda/2,\,
Nf(x)\leq\delta\lambda\}.
\end{equation}

Let $k\in \Z$ be such that $Q_i\in \AD_{k}$.
Let us check that
\begin{equation}  \label{claim1}
S_{k+3}\chi_{A_\lambda}(y) \leq C\,\delta\ve^{-1}\qquad
\mbox{for all $y\in Q_i$}.
\end{equation}
For a fixed $y\in Q_i$, let $j_0$ be the least non negative integer such that
there exists some $y_0\in 2^{j_0}Q_{y,k+3}\cap A_\lambda$. Let us denote $C_j =
2^jQ_{y,k+3} \setminus 2^{j-1}Q_{y,k+3}$ for $j>j_0$, and
$C_{j_0} =2^{j_0}Q_{y,k+3}$. Then we have
$$S_{k+3}\chi_{A_\lambda}(y) = \int_{A_\lambda\cap 2Q_i} s_{k+3}(y,z)\,
d\mu(z)
\leq C\,\sum_{j=j_0}^{n_0} \frac{\mu(A_\lambda\cap C_j)}{
\ell(2^j Q_{y,k+3})^n},$$
where $n_0$ is the least integer such that $2Q_{y,k+2}\subset 2^{n_0}
Q_{y,k+3}$.
Let $V_j$ be the $\ell(2^j Q_{y,k+3})$-neighborhood of $C_j$. It is easily
checked that $T_* (f\chi_{3Q_i\setminus V_j})(z) \leq C\,Nf(z)$ for all
$z\in C_j$. Therefore, if $\delta$ is small enough, for $z\in A_\lambda\cap
C_j$ we must have $T_*(f\chi_{V_j})(z)\geq \ve\lambda/4$.
Then, by the weak (1,1) boundedness of $T_*$ with respect to $\mu$, we get
$$\mu(A_\lambda\cap C_j) \leq \mu\{z:\,T_*(f\chi_{V_j})(z)\geq \ve\lambda/4\}
\leq \frac{C}{\ve \lambda} \int_{V_j} |f|\,d\mu.$$
Using the finite overlap of the neighborhoods $V_j$,
$$S_{k+3}\chi_{A_\lambda}(y) \leq  \frac{C}{\ve \lambda}
\sum_{j=j_0}^{n_0} \frac{1}{\ell(2^j Q_{y,k+3})^n} \int_{V_j} |f|\,d\mu
\leq \frac{C}{\ve \lambda} Nf(y_0) \leq \frac{C\delta }{\ve},$$
which proves \rf{claim1}.

By \rf{mmhh2}, we get
$S_{k+3} \chi_{\R^d\setminus A_\lambda} (y) \geq 9/10 - C\delta\ve^{-1}
>1/4$ for all $y\in Q_i$, if $\delta$ is small enough.
By the $Z_\infty$ condition,
$w(2Q_i\setminus A_\lambda) \geq \tau w(Q_i)$.
Therefore, by the finite overlap of the cubes $2Q_i$,
\begin{equation} \label{nn*2}
w(\Omega_\lambda) \leq \tau^{-1} \sum_i w(2Q_i\setminus A_\lambda)
\leq C\,\alpha^{-1} w(\Omega_\lambda\setminus A_\lambda).
\end{equation}
Thus,
$w(A_\lambda) \leq (1-C\,\tau^{-1})\,w(\Omega_\lambda)$.
Now we only have to take $\eta :=  1-C\,\tau^{-1}$ (which does not depend
on $\delta$, $\ve$ or $\lambda$), and \rf{gli1} follows.
\end{proof}

% ***************************************************************************
% ***************************************************************************
% ***************************************************************************

\section{The implications {\rm(b)} $\Rightarrow$ {\rm(a)} $\Rightarrow$
 {\rm (d)} $\Rightarrow$  {\rm (e)} of Lemma \ref{mainw}}

The implication  {\rm(b)} $\Rightarrow$ {\rm(a)} is trivial. Let us
see the remaining ones.

\begin{proof}[{\bf Proof of (a) $\Rightarrow$ (d) in Lemma \ref{mainw}}]
We have defined the kernels $s_k(x,y)$ so that they are CZ kernels
uniformly on $k\in\Z$. By the statement (a) in Lemma \ref{mainw} we
know that they are of weak type $(p,p)$ with respect to $w\,d\mu$.
We only have to check that this holds {\em uniformly} on $k$.
Indeed, if this is not the case, for each $m\geq1$
we take $S_{k_m}$ such that $\|S_{k_m}\|_{L^p(w),L^{p,\infty}(w)} \geq
m^3$. Then we define $T=\sum_{m=1}^\infty \frac{1}{m^2}S_{k_m}$. Since
$\sum_m  \frac{1}{m^2}<\infty$, $T$ is a CZO (using also uniform
estimates for the operators $S_k$). On the other hand, we have
$\|T\|_{L^p(w),L^{p,\infty}(w)}\geq\frac{1}{m^2}
\|S_{k_m}\|_{L^p(w),L^{p,\infty}(w)}\geq m$ for each $m$, because $S_k$ are
integral operators with non negative kernel. Thus
$\|T\|_{L^p(w),L^{p,\infty}(w)}=\infty$, which contradicts the
statement (a) in Lemma \ref{mainw}.
\end{proof}

\vspace{2mm}
\begin{proof}[{\bf Proof of (d) $\Rightarrow$ (e) in Lemma
\ref{mainw} for $1<p<\infty$.}]
Since the operators $S_k$ are of weak type $(p,p)$  with respect to $w\,d\mu$,
from \rf{quho} it follows that their duals are also of weak type
$(p,p)$ with respect to $w\,d\mu$, uniformly on $k$.
Then, the statement (e) is a consequence
of duality in Lorentz spaces. We only have to argue as in
\cite[p.\ 341]{Sawyer0.5}, for example:
\begin{eqnarray*}
\left(\int |S_k(w\,\chi_Q)|^{p'}\,\sigma\,d\mu \right)^{1/p'}\!\! & = &
\sup_{\|f\|_{L^p(\sigma)}\leq1} \int S_k(w\,\chi_Q)\,f\,\sigma\,d\mu\\
& = & \sup_{\|f\|_{L^p(\sigma)}\leq1} \int_Q S_k^*(f\,\sigma)\,w\,d\mu \\
& = & \sup_{\|f\|_{L^p(\sigma)}\leq1} \int_0^\infty \!
w\{x\in Q:\,S_k^*(f\,\sigma)(x)>\lambda\}\, d\lambda \\
& \leq & \int_0^\infty \min\bigl(C\,\lambda^{-p},\, w(Q)\bigr)\,d\lambda \,=\,
C\,w(Q)^{1/p'}.
\end{eqnarray*}
\end{proof}

% ***************************************************************************
% ***************************************************************************
% ***************************************************************************

\section{The implication {\rm(e)} $\Rightarrow$ {\rm (c)} of Lemma
\ref{mainw}}

We need to introduce some notation and terminology. Let $\Omega$ be an open
set. Suppose that we have a Whitney decomposition $\Omega=\bigcup_i Q_i$
into dyadic cubes $Q_i$ with disjoint interiors, with $10Q_i\subset \Omega$,
$\dist(Q_i,\partial\Omega)\approx
\ell(Q_i)$, and such that the cubes $4Q_i$ have finite overlap.
We say that two cubes $Q$ and $R$ are neighbors if $Q\cap R\neq
\varnothing$ (recall that we are assuming that the cubes are closed).
For a fixed $i$, we denote by $U_1(Q_i)$ the union of all the neighbors
of $3Q_i$ (including $Q_i$ itself). For $m>1$, inductively we let
$U_m(Q_i)$ be the union of all the cubes which are neighbors of some cube in
$U_{m-1}(Q_i)$. That is, one should think that $U_m(Q_i)$ is formed by $3Q_i$
and $m$ ``layers'' of neighbors.

We denote by $M_R$ the non centered radial maximal
Hardy-Littlewood operator:
$$M_R f(x) = \sup_{B:\,x\in B} \frac{1}{r(B)^n} \int_B |f|\,d\mu,$$
where $B$ stands for any ball containing $x$ and $r(B)$ is its radius.

In order to prove the implication {\rm(e)} $\Rightarrow$ {\rm (c)} we will
need a very sharp maximum principle. In the following lemma we deal with this
question.

\begin{lemma}  \label{prmax}
Let $\ve>0$ be some arbitrary fixed constant. There exist $\beta>0$ and
$m\geq1$, $m\in\Z$, both big enough, such that the operator
$T = N + \beta\, M_R$
satisfies the following maximum principle for all $\lambda>0$ and all
$f\in L^1_{loc}(\mu)$:
Let $\Omega_\lambda = \{x:\, Tf(x)>\lambda\}$, and consider a Whitney
decomposition $\Omega_\lambda=\bigcup_i Q_i$
as above. Then, for any $x\in Q_i$,
\begin{equation}  \label{eqprmax}
T(f\,\chi_{\R^d\setminus U_m(Q_i)})(x) \leq (1+\ve)\,\lambda.
\end{equation}
\end{lemma}

The point in this lemma is that the constant $\ve>0$ can be as small as we
need, which will be very useful.
We only have to define the operator $T$ choosing $\beta$ big enough, and also
to take the integer $m$ sufficiently big in $U_m(Q_i)$.
Notice also that
$Nf(x) \leq Tf(x) \leq (1+C\,\beta) Nf(x)$,
because $M_R f(x) \leq C\,Nf(x)$.

\vspace{3mm}
\begin{proof}
Let $x\in Q_i$ be some fixed point.
First we will show that, for some $z\in\partial\Omega$,
\begin{equation} \label{max1}
M_R (f\,\chi_{\R^d\setminus U_m(Q_i)})(x) \leq (1+\ve/2)\,M_R f(z),
\end{equation}
if we choose $m$ big enough.
Let $B$ be some ball containing $x$ such that
$$(1+\ve/2)^{1/2}\,
\frac{1}{r(B)^n}\int_B |f\,\chi_{\R^d\setminus U_m(Q_i)}|\,d\mu
\geq M_R (f\,\chi_{\R^d\setminus U_m(Q_i)})(x).$$
Notice that if $M_R (f\,\chi_{\R^d\setminus U_m(Q_i)})(x)\neq0$,
then $B\setminus U_m(Q_i) \neq \varnothing$.
Since $3Q_i\subset U_m(Q_i)$, we get
\begin{equation} \label{inf0}
\diam(B) \geq \ell(Q_i).
\end{equation}
Recall also that $U_m(Q_i)$ is formed by $m$ ``layers'' of Whitney cubes,
and so we have
\begin{equation} \label{infm}
\diam(B) \geq m\!\!\inf_{
\mbox{\scriptsize $
\begin{array}{c} j:\,Q_j\subset U_m(Q_i)\\
Q_j\cap B\neq \varnothing\end{array}$}} \!\! \ell(Q_j).
\end{equation}

We distinguish two cases:

\begin{itemize}

\item[(a)] Assume $100\,\ell(Q_i) \leq [(1+\ve/2)^{1/2n} -1]\,r(B) =:
 C_\ve \,r(B)$. That is,
$\ell(Q_i)$ is small compared to $r(B)$.
We choose $z\in\partial
\Omega$ such that $\dist(x,\R^d\setminus\Omega) = |x-z| \leq 100\ell(Q_i)$.
Then there exists some ball $B'$ containing $z$ and $B$ with radius
$r(B') \leq r(B) + |x-z| \leq (1+\ve/2)^{1/2n}\,r(B)$.
Therefore,
\begin{equation} \label{casa}
M_Rf(z)\geq \frac{1}{r(B')^n} \int_{B'} |f|\,d\mu \geq
\frac{1}{(1+\ve/2)^{1/2}\,r(B)^n} \int_{B} |f|\,d\mu,
\end{equation}
and \rf{max1} holds.

\item[(b)] Suppose that
$100\,\ell(Q_i) \geq C_\ve \,r(B)$.
Then there exists some Whitney cube $P$ in $U_m(Q_i)$ such that $P\cap B\neq
\varnothing$ and
$\ell(P) \leq 300\, C_\ve^{-1}\ell(Q_i)/m$. Otherwise, by \rf{infm},
$2\,r(B) \geq 300 \,C_\ve^{-1}\ell(Q_i),$
which contradicts our assumption.

Since $P\cap B\neq \varnothing$,
we can find $z\in\partial\Omega$ such that $\dist(z,B)\leq 100\,\ell(P)$. Thus,
$$\dist(z,B)\leq  30000 \,C_\ve^{-1}\,\ell(Q_i)/m \leq C_\ve\,\ell(Q_i)/2,$$
if $m$ is chosen big enough. By \rf{inf0}, we obtain
$\dist(z,B)\leq C_\ve\,r(B)$. Then there exists some ball $B'$ containing
$z$ and $B$ with radius
$$r(B') \leq (1+C_\ve)\,r(B) = (1+\ve/2)^{1/2n}\,r(B).$$
Arguing as in \rf{casa}, we obtain \rf{max1}.
\end{itemize}

Now we have to deal with the term $Nf(x)$.
Notice that if $z\in\partial \Omega$ is the point
chosen in (a) or (b) above, then in both cases we have
$|x-z|\leq C\,\ell(Q_i)$, where $C$ may depend on $m$.
Thus we may choose some constant  $\eta>0$ big enough so that
$\eta\,\ell(Q_i)\gg \dist(x,\partial\Omega)$, $|x-z|$. We set
$B_\eta:= B(x,\eta\,\ell(Q_i))$, and we have
$$N (f\,\chi_{\R^d\setminus U_m(Q_i)})(x) \leq
N (f\,\chi_{B_\eta\setminus U_m(Q_i)})(x) +
N (f\,\chi_{\R^d\setminus B_\eta})(x).$$
Since $|x-z|\ll \eta \,\ell(Q_i)$, for each $k$ we get
$$|S_k(f\,\chi_{\R^d\setminus B_\eta})(x) - S_k(f\,\chi_{\R^d\setminus
B_\eta})(z)| \leq C_{12}\,M_Rf(z),$$
where $C_{12}$ may depend on $\eta$. Thus
$N(f\,\chi_{\R^d\setminus B_\eta})(x) \leq
Nf(z) + C_{12}\,M_Rf(z)$.
We also have
$N (f\,\chi_{B_\eta\setminus U_m(Q_i)})(x) \leq C_{13}  \,M_R f(z),$
with $C_{13}$ depending on $\eta$.
Therefore,
$$N (f\,\chi_{\R^d\setminus U_m(Q_i)})(x)
\leq N f(z) + C_\eta  \,M_R f(z).$$
If we take $\beta$ such that $C_\eta\leq  \beta\ve/2$,
by \rf{max1}, we obtain
\begin{eqnarray*}
T(f\,\chi_{\R^d\setminus U_m(Q_i)})(x) & \leq &
N f(z) + C_\eta  \,M_R f(z) + \beta\,(1+\ve/2)\,M_R f(z) \\
& \leq & N f(z) + \beta\,(1+\ve)\,M_R f(z) \\
& \leq & (1+\ve)\,Tf(z) \,\leq \,(1+\ve)\,\lambda.
\end{eqnarray*}
\end{proof}

% ****************************************************************************
% ****************************************************************************

\vspace{3mm}
\begin{proof}[{\bf Proof of (e) $\Rightarrow$ (c) in Lemma
    \ref{mainw} for $1<p<\infty$.}]
We will show that for some $\beta\geq0$, the operator
$T:= N + \beta\,M_R$ is bounded on $L^p(w)$.
The precise value of $\beta$ will be fixed below.
Without loss of generality, we take $f\in L^1(\mu)$ non negative
with compact support.
Given any $\lambda>0$, we denote
$\Omega_\lambda = \{x:\, Tf(x)>\lambda\}.$
We will show that there exists some constant $\eta$, with
$0<\eta<1$, such that for all $\ve,\lambda>0$
\begin{equation}\label{comgli}
w(\Omega_{(1+\ve)\lambda}) \leq \eta\,w(\Omega_\lambda) +
\frac{C_\ve}{\lambda^p} \int |f|^p\,w d\mu,
\end{equation}
where $C_\ve$ is some constant depending on $\ve$ but not on $\lambda$.
It is straightforward to check that
\rf{comgli}
implies that $T$ is of weak type $(p,p)$ with respect to $w\,d\mu$ for
$\ve$ small enough.

As in Lemma \ref{prmax}, we consider the Whitney decomposition
$\Omega_\lambda = \bigcup_i Q_i$, where $Q_i$ are dyadic cubes with disjoint
interiors (the {\it Whitney cubes}).
Suppose that $m$ and $\beta$ are chosen in Lemma \ref{prmax} so that the
maximum principle \rf{eqprmax} holds with $\ve/2$ instead of $\ve$.
Take some cube $Q_i\subset \Omega_\lambda$.
To simplify notation, we will write $U_i$ instead of $U_m(Q_i)$.
Then, for $x\in Q_i\cap \Omega_{(1+\ve)\lambda}$, we have
$T(f\,\chi_{\R^d\setminus U_i})(x) \leq (1+\ve/2)\,\lambda$,
and so
\begin{equation} \label{zzz2}
T(f\,\chi_{U_i})(x) \geq \ve \, \lambda/2.
\end{equation}
Let $h\in \Z$ be such that $Q_i\in \AD_{h}$.
If for all $k$ with $h-n_1\leq k\leq h+5$ we have
$S_k (f\chi_{U_i}) (x) \leq \delta\,\lambda$,
where $n_1,\delta$ are positive constants which we will fix below, then
we write $x\in B_\lambda$ (i.e.\ $x$ is a ``bad point'') and, otherwise, $x\in
G_\lambda$.

Notice that $G_\lambda \cup B_\lambda = \Omega_{(1+\ve)\lambda} \subset
\Omega_\lambda$. We will see that $B_\lambda$ is quite small.
Indeed, we will prove that
\begin{equation} \label{blb1}
w(B_\lambda) \leq \eta_1 w(\Omega_\lambda),
\end{equation}
for some positive constant $\eta_1<1$ independent of $\ve$ and $\lambda$.

Assume that \rf{blb1} holds for the moment, and let us estimate $w(G_\lambda)$.
For $Q_i\in \AD_{h}$, we have
\begin{eqnarray} \label{eqhy1}
w(Q_i \cap G_\lambda) & \leq & \frac{1}{\delta\lambda} \int_{Q_i}
\sum_{k=h-n_1}^{h+5} S_k(f\,\chi_{U_i})\,w\,d\mu \nonumber \\
& = & \frac{1}{\delta\lambda} \sum_{k=h-n_1}^{h+5}
\int f\,\chi_{U_i} S_k(w\,\chi_{Q_i})\,d\mu \nonumber \\
& \leq &  \frac{1}{\delta\lambda} \sum_{k=h-n_1}^{h+5}
\left(\int S_k(w\,\chi_{Q_i})^{p'}\,\sigma\,d\mu\right)^{1/p'}
\left(\int_{U_i} |f|^p\,w\, d\mu\right)^{1/p}\nonumber \\
& \leq & \frac{C(n_1+6)}{\lambda} w(Q_i)^{1/p'}
\left(\int_{U_i} |f|^p\,w\, d\mu\right)^{1/p}.
\end{eqnarray}
Using the inequality $a^{1/p'}\,b^{1/p}\leq \theta a + \theta^{-p'/p} b$, for
$a,\,b,\,\theta>0$, we get
$$w(Q_i \cap G_\lambda) \leq \theta\, w(Q_i) +
\frac{C\theta^{-p'/p}}{\lambda^p} \int_{U_i} |f|^p\,w\, d\mu.$$
It is not difficult to check that the family of sets $\{U_i\}_i$ has
bounded overlap (depending on $m$). Then,
summing over all the indices $i$, we obtain
$$w(G_\lambda) \leq C\theta\, w(\Omega_\lambda) + \frac{C(\theta,m)}{\lambda^p}
\int |f|^p\,w\, d\mu.$$
By \rf{blb1}, if we choose $\theta = (1-\eta_1)/2C$, we get
$$w(\Omega_{(1+\ve)\lambda}) \leq \frac{1+\eta_1}{2} \,w(\Omega_\lambda) +
\frac{C}{\lambda^p} \int |f|^p\,w\, d\mu,$$
which is \rf{comgli} with $\eta=(1+\eta_1)/2$.

Now we have to show that \rf{blb1} holds.
We intend to use the $Z_\infty$ property. Let us see that
\begin{equation} \label{zzz1}
S_{h+3} \chi_{\R^d\setminus B_\lambda}(y) \geq \frac{1}{4}
\end{equation}
for all $y\in Q_i$. By \rf{zzz2}, if $z\in Q_i$, then
$N(f\,\chi_{U_i})(z) \geq C_{14} \, \lambda$, where $C_{14}$ is some
positive constant depending on $\ve,\,\beta$. Then we have
\begin{equation} \label{nn*10}
S_k (f\,\chi_{U_i})(z) \geq C_{14} \, \lambda
\end{equation}
for some $k\geq h-n_1$. If moreover $z\in B_\lambda\cap Q_i$, then
this inequality holds for some $k\geq h+6$, assuming
that we take $\delta<C_{14}$.

Suppose that $B_\lambda\cap \supp(s_{h+3}(y,\cdot))\neq\varnothing$.
Let $j_0\geq0$ be the least
integer such there exists some $x_0\in 2^{j_0}Q_{y,h+3}$, and let $n_0$ be
the least integer such that $Q_{y,h+2}\subset 2^{n_0} Q_{y,h+3}$.
Then we have
\begin{eqnarray*}
S_{h+3} \chi_{B_\lambda}(y) & = &
\int_{z\in B_\lambda} s_{h+3}(y,z)\,d\mu(z) \\
& \leq & C\, \sum_{j=j_0}^{n_0}
\frac{\mu(B_\lambda\cap (2^{j+1}Q_{y,h+3} \setminus 2^{j}Q_{y,h+3})
)}{\ell(2^jQ_{y,h+3})^n}.
\end{eqnarray*}
It is not difficult to check that if $z\in B_\lambda\cap (2^{j+1}Q_{y,h+3}
\setminus 2^{j}Q_{y,h+3})$ and $k\geq h+6$, then $\supp(s_k(z,\cdot))\subset
2^{j+2}Q_{y,h+3} \setminus 2^{j-1}Q_{y,h+3} =: V_j$. Therefore,
$N(f\chi_{V_j})(z) \geq  C_{14} \, \lambda$. Then,
by the weak $(1,1)$ boundedness of $N$, we have
\begin{eqnarray*}
\mu(B_\lambda\cap (2^{j+1}Q_{y,h+3} \setminus 2^{j}Q_{y,h+3})) & \leq &
\mu\{z:\, N(f\,\chi_{V_j})(z)\geq C_{14}\lambda\}\\
& \leq & \frac{C}{\lambda} \int_{V_j}|f|\,d\mu.
\end{eqnarray*}
Thus, by the finite overlap of the sets $V_j$, and since $x_0\in B_\lambda$,
\begin{eqnarray*}
S_{h+3} \chi_{B_\lambda}(y) & \leq & \frac{C}{\lambda}\,
\sum_{j=j_0}^{n_0}
\frac{1}{\ell(2^jQ_{y,h+3})^n} \int_{V_j}|f|\,d\mu  \\
& \leq & \frac{C}{\lambda}
\,(S_{h+2} f(x_0) + S_{h+3} f(x_0)+ S_{h+4} f(x_0)) \, \leq\, C_{15}\delta.
\end{eqnarray*}
Notice that $C_{15}$ depends on $\ve$, but not on $\delta$.
If $\delta$ is small enough, we obtain $S_{h+3} \chi_{B_\lambda}(y) \leq
1/4$. Now, we have
$S_{h+3} \chi_{\R^d\setminus B_\lambda}(y) \geq 9/10
- S_{h+3} \chi_{B_\lambda}(y) \geq 1/4$,
and \rf{zzz1} holds.

By the $Z_\infty$ property,  we get
$w(2Q_i\setminus B_\lambda) \geq \tau\, w(Q_i)$,
and because of the finite overlap of the cubes $2Q_i$,
$$w(\Omega_\lambda)  \leq \tau^{-1}
 \sum_i w(2Q_i\setminus B_\lambda) \leq C_{16} \tau^{-1}
 w(\Omega_\lambda\setminus B_\lambda),$$
which implies \rf{blb1}.
\end{proof}

A slight modification of the arguments above yields the
{\bf proof of the implication (d) $\Rightarrow$ (c) in the weak
$(1,1)$ case}. Instead of using \rf{cond1'} to estimate $w(Q_i \cap
G_\lambda)$ in \rf{eqhy1}, one can apply directly that the operators $S_k$ are
bounded from $L^1(w)$ into $L^{1,\infty}(w)$. We leave the details for
the reader.

% ***************************************************************************
% ***************************************************************************

% ***************************************************************************
% ***************************************************************************

% ***************************************************************************
% ***************************************************************************

\section{Preliminary lemmas for the proof of Lemma \ref{main}}

Section 8--10 are devoted to the proof of Lemma \ref{main}.
As in our lemma about the weak $(p,p)$ case, the
implication (c) $\Rightarrow$ (b) is a direct consequence of the good
$\lambda$ inequality of Lemma \ref{lemagli}. The proofs of the implications
(b) $\Rightarrow$ (a) $\Rightarrow$ (d) $\Rightarrow$ (e) are similar
to the ones of Lemma \ref{mainw}. We will not go through the details.
So we have to concentrate on the implication (e) $\Rightarrow$ (c),
which is more difficult than the corresponding implication of the weak
$(p,p)$ case, as we will see.

In this section we will obtain some technical results which will be
needed later for the proof of (e) $\Rightarrow$ (c).

% ****************************************************************************

\begin{lemma} \label{perprmax}
Let $\rho\geq1$ be some fixed constant. Let $Q\subset\R^d$ be some cube and
suppose that $x\in Q\cap\supp(\mu)$, $x'\in \rho Q\cap\supp(\mu)$, and
$y\in\R^d\setminus 2Q$. Then, $s_k(x,y)\leq C_{17}\sum_{j=k-5}^{k+5}
s_j(x',y)$,
for any $k\in\Z$, with $C_{17}$ depending on $\rho$ and
assuming $A$ big enough (depending on $\rho$ too).
\end{lemma}

\begin{proof}
We denote $s^k(x',y)= \sum_{j=k-5}^{k+5} s_j(x',y)$.
Observe that, by the definition of $s_j(x',y)$, we have
\begin{equation} \label{ssupk}
s^k(x',y) \approx \min\left(\frac1{A\ell(Q_{x',k+5})^n},\,
\frac1{A|x'-y|^n}\right) \qquad \mbox{if $y\in Q_{x',k-5}$}.
\end{equation}

Let $h\in\Z$ be such that $Q\in\AD_h$. If $k\geq h+3$, then
$\supp(s_k(x,\cdot))\subset 2Q$ by Lemma \ref{suportpetit}, and then
$s_k(x,y)=0$.

Assume now $k\leq h-3$. Since $Q\in\AD_h$, we have $Q\subset Q_{x,h-1}$.
If $A$ is big enough (depending on $\rho$), we deduce
$x'\in Q_{x,h-2}\subset
Q_{x,k}$ by (g) in Lemma \ref{lattice}.
Then we get $2Q_{x,k-1}\subset Q_{x',k-4}$, and so $y\in
Q_{x',k-4}$ if $s_k(x,y)\neq0$.
We also deduce
$\ell(Q_{x',k+5})\ll\ell(Q_{x,k})$. By \rf{ssupk}, if $s_k(x,y)\neq0$,
we obtain
\begin{eqnarray*}
s^k(x',y)  &\geq & C^{-1}
\min\left(\frac1{A\ell(Q_{x,k})^n},\,
\frac1{A|x'-y|^n}\right) \\
&\geq & C^{-1}
\min\left(\frac1{A\ell(Q_{x,k})^n},\,
\frac1{A|x-y|^n}\right) \, \geq \,s_k(x,y).
\end{eqnarray*}

Suppose finally that $|h-k|\leq2$. As above, we have $x'\in Q_{x,h-2}$, and
since $h-2\geq k-4$, $x'\in Q_{x,k-4}$.
Then we get $2Q_{x,k-1}\subset Q_{x',k-5}$, and so $y\in
Q_{x',k-5}$ if $s_k(x,y)\neq0$.
On the other hand, $Q\not\subset Q_{x',h+1}$, because $Q\in\AD_h$. Thus
$$\ell(Q_{x',h+1}) \leq C\ell(Q) \leq C|x-y|,$$
with $C$ depending on $\rho$. Then, if $s_k(x,y)\neq0$, by \rf{ssupk} we get
$$s^k(x',y) \geq C^{-1} \frac1{A|x-y|^n} \geq C^{-1}\,s_k(x,y).$$
\end{proof}

% ****************************************************************************

Given $\alpha,\beta>1$,
we say that some cube $Q\subset\R^d$ is
$\mu$-$\sigma$-$(\alpha,\beta)$-doubling if $\mu(\alpha Q)\leq \beta\,
\mu(Q)$ and  $\sigma(\alpha Q)\leq \beta\,\sigma(Q)$. We say that $Q$
is $\mu$-$\sigma$-doubling if $\alpha=2$ and $\beta$ is some fixed
constant big enough (which perhaps is not specified explicitly).
Next lemma deals with the existence of this kind of cubes.

\begin{lemma}  \label{dobsim}
Suppose that the operators $S_k$ are bounded on $L^r(\sigma)$ uniformly on $k$
for some $r$ with $1<r<\infty$ and that the constant $A$ is big enough.
Then there exists some constant $\beta>0$ such that
for any $x\in\supp(\mu)$ and
all cubes $Q,R$ centered at $x$ with $\delta(Q,R)\geq A/2$, there exists
some $\mu$-$\sigma$-$(100,\beta)$-doubling cube $P$ centered at $x$,
with $Q\subset P\subset R$.
\end{lemma}

\begin{proof}
The constant $\beta$ will be chosen below. For the moment, let
us say that $\beta\geq 100^{d+1}$.
Let $N_0$ be the least integer such that
$R\subset 100^{N_0}Q$. For each $j\geq0$,
we denote $R_j=100^{-j}R$. We have $\delta(R_{N_0},R) \geq A/2 -
100C_0 > A/4$. We will show that some cube $R_j$, with $j=0,\ldots
N_0$, is doubling with respect to $\mu$ and $\sigma$.

Suppose that all the cubes $R_j$, $j=0,\ldots
N_0$, are either non $(100,\beta)$-$\mu$-doubling, or non
$(100,\beta^{r})$-$\sigma$-doubling (for simplicity, we will show the
existence of a cube $P$ which is $(100,\beta^r)$-$\sigma$-doubling, instead of
$(100,\beta)$-$\sigma$-doubling).
For each $j=0,\ldots,N_0$, let $a_j$ be the number of non
$(100,\beta)$-$\mu$-doubling cubes of the form $100^{-k}R$,
$k=0,\ldots,j$ and let $b_j$ the number of
non $(100,\beta^{r})$-$\sigma$-doubling cubes of the form
$100^{-k}R$, $k=0,\ldots,j$.
From our assumption we deduce
\begin{equation}  \label{ajbj}
a_j + b_j \geq j+1.
\end{equation}
We have $\mu(R_j) \leq \beta^{-a_j} \mu(R)$. Thus,
\begin{equation} \label{nndd1}
\frac{\mu(R_j)}{\ell(R_j)^n} \leq \frac{\beta^{-a_j}\mu(R)}{
100^{-jn}\ell(R)^n} \leq C_{0}\, \frac{100^{jn}}{\beta^{a_j}}.
\end{equation}

Let $R_{s_0}$ be the largest non $(100,\beta^{r})$-$\sigma$-doubling
cube of the form
$100^{-k}R$, $k=0,\ldots,N_0$. Then, for $j\geq
s_0$ we have
$$\sigma(R_j)\leq  \beta^{-rb_j} \sigma(100R_{s_0}) \leq \frac{1}{2}\,
\beta^{-rb_j} \sigma(100R_{s_0}\setminus R_{s_0}),$$
if $\beta$ is big enough.

Let $h\in\Z$ be such that $Q\in \AD_{h}$. We denote
$S= \sum_{i=-3}^3 S_{h+i}$.
From the properties of the
kernels $s_k(x,y)$, it is easily seen that, for $x\in
100R_{s_0}\setminus R_{s_0}$ and $j=s_0,s_0+1,\ldots,N_0$, we have
$S(\chi_{R_j})(x) \geq C^{-1}\mu(R_j)/\ell(R)^{n}$.
Then, using the boundedness of $S$ on $L^r(w)$, we obtain
\begin{eqnarray*}
C\,\sigma(R_j) \,\geq\,
\int_{100R_{s_0}\setminus R_{s_0}} |S(\chi_{R_j})|^r\,\sigma\,d\mu & \geq &
C^{-1}\, \sigma(100R_{s_0}\setminus R_{s_0})\,
\frac{\mu(R_j)^r}{\ell(R)^{nr}} \\
& \geq & C^{-1} \, \frac{\beta^{b_jr}}{100^{jnr}}\, \sigma(R_j)\,
 \frac{\mu(R_j)^r}{\ell(R_j)^{nr}},
\end{eqnarray*}
if $j\geq s_0$. Thus,
\begin{equation} \label{nndd2}
\frac{\mu(R_j)}{\ell(R_j)^{n}} \leq C\,
\frac{100^{jn}}{\beta^{b_j}}\qquad \mbox{if $j\geq s_0$.}
\end{equation}

By \rf{ajbj}, $\max(a_j,\,b_j)\geq (j+1)/2$. Then,
from \rf{nndd1} and \rf{nndd2}, we get
$$\frac{\mu(R_j)}{\ell(R_j)^{n}} \leq C\,
\frac{100^{jn}}{\beta^{(j+1)/2}}.$$
Therefore,
$$\delta(R_{N_0},R) \leq \sum_{j=0}^\infty
C\,\left(\frac{100^{n}}{\beta^{1/2}}\right)^{j+1} \leq
\frac{C}{\beta^{1/2}},$$
if $\beta^{1/2}>2\!\cdot\!100^n$.
Thus $\delta(R_{N_0},R) \leq A/4$ if $\beta$ is big
enough, which is a contradiction.
\end{proof}

Let us remark that if in the lemma above we also assume that the
operators $S_k$ are bounded uniformly on $k$ on $L^s(w)$ for some $s$ with
$1<s<\infty$, then it is possible to show the existence
of cubes which are $\mu$-doubling, $\sigma$-doubling and $w$-doubling
simultaneously, by an easy modification of the proof.

Notice also that if $\int|S_k(\sigma\chi_Q)|^p\,w\,d\mu\leq C\,\sigma(Q)$
for $k\in\Z$ and all the cubes $Q\subset\R^d$, then $N$
is of weak type $(p',p')$
with respect to $\sigma$ and bounded on $L^r(\sigma)$ for
$p'<r\leq\infty$. Thus the assumptions of the preceding lemma are satisfied.

% ****************************************************************************

\begin{lemma}  \label{sigmenor}
Suppose that the operators $S_k$ are bounded on $L^r(\sigma)$ uniformly on $k$
for some $r$ with $1<r<\infty$ and that the constant $A$ is big enough.
Then there exists some constant $\eta$ with $0<\eta<1$ such that,
for all $x\in \R^d$ and $k\in\Z$, $\sigma(Q_{x,k}) \leq \eta
\sigma(Q_{x,k-1})$.
\end{lemma}

\begin{proof}
We denote $S=\sum_{i=-2}^2 S_{h+i}$. Then, we have
$S(\chi_{Q_{x,k-1}\setminus Q_{x,k}})(y) \geq C_{18}>0$,
for all $y\in Q_{x,k}$. Therefore,
$$\sigma(Q_{x,k}) \leq C_{18}^{-r} \int |S(\chi_{Q_{x,k-1}\setminus
Q_{x,k}})|^r\,\sigma\,d\mu \leq C_{19}\,\sigma(Q_{x,k-1}\setminus Q_{x,k}).$$
Thus, $\sigma(Q_{x,k-1}) \geq (1+C_{19}^{-1})\,\sigma(Q_{x,k})$.
\end{proof}

% ****************************************************************************

We will use the following version of Wiener's Covering Lemma.

\begin{lemma}\label{vitali}
Let $A\subset\R^d$ be a bounded set and $\{Q_i\}_{i\in I}$
some family of cubes such that $A\subset \bigcup_{i\in I}Q_i$,
with $Q_i\cap A\neq \varnothing$ for each $i\in I$.
Then there exists some finite or countable subfamily $\{Q_j\}_{j\in J}$,
$J\subset I$, such that
\begin{itemize}
\item[(1)] $A\subset \bigcup_{j\in J} 20 Q_j$,
\item[(2)] $2Q_j\cap 2Q_k = \varnothing$ if $j,k\in J$.
\item[(3)] If $j\in J$, $k\not\in J$, and $2Q_j\cap 2Q_k\neq\varnothing$,
then $\ell(Q_k)\leq 10\ell(Q_j)$.
\end{itemize}
\end{lemma}

The main novelty with respect to the usual Wiener's
Lemma is the assertion (3).
Although the lemma follows by standard arguments, for the sake
of completeness we will show the detailed proof.

\begin{proof}
We will construct inductively the set $J=\{j_1,j_2,\ldots\}$.
Let $\ell_1=\sup_{i\in I} \ell(Q_i)$. If $\ell_1=\infty$,
the lemma is straightforward. Otherwise,
we take $Q_{j_1}$ such that $\ell(Q_{j_1})>d_1/2$.
Assume that $Q_{j_1},\ldots,Q_{j_{m-1}}$ have been chosen. We set
$$\ell_{m}
= \sup\Bigl\{\ell(Q_i):\,4Q_i\not\subset{\ts \bigcup_{s=1}^{m-1}}
20Q_{j_s}\Bigr\},$$
and we choose $Q_{j_{m}}$ such that $\ell(Q_{m})>\ell_{m}/2$ and
$4Q_{j_{m}}\not\subset\bigcup_{s=1}^{m-1} 20Q_{j_s}$.

By construction, $A\subset \bigcup_{m=1}^\infty 20Q_{j_m}$, and also
$\ell(Q_{j_{m}})\geq \ell(Q_{j_s})/2$ for $s>m$.
We have $2Q_{j_m}\cap 2Q_{j_s}=\varnothing$ for each $s=1,\ldots,m-1$,
because otherwise $2Q_{j_m}\subset 10Q_{j_s}$, and then
$4Q_{j_m}\subset 20Q_{j_s}$.

Finally we show that the third property holds. Suppose that
$k\not\in J$ and $2Q_{j_m}\cap 2Q_k\neq\varnothing$.
If $\ell(Q_k)> 10\ell(Q_{j_m})$, it is easily seen that
$4Q_{j_m}\subset4Q_k$. Because of the definition of $Q_{j_m}$, we must have
$4Q_k\subset \bigcup_{s=1}^{m-1}20Q_{j_s}$ (otherwise $\ell_m\geq
\ell(Q_k) > 10\ell(Q_{j_m})$,  which is not possible). However the last
inclusions imply
$4Q_{j_m}\subset \bigcup_{s=1}^{m-1}20Q_{j_s}$, which is a contradiction.
\end{proof}

% ****************************************************************************
% ****************************************************************************
% ****************************************************************************

\section{Boundedness of $N$ over functions of type  $\sigma\chi_Q$ on
$L^p(w)$}

The main result of this section is the following lemma.

% ****************************************************************************

\begin{lemma}\label{fonam**}
If
$$\int |S_k(\sigma\,\chi_Q)|^p\, w\,d\mu \leq C\,\sigma(Q)$$
for all cubes $Q\subset \R^d$ uniformly on $k\in\Z$, then
$$\int |N(\sigma\,\chi_Q)|^p\, w\,d\mu \leq
C\,\sigma({\ts\frac{11}{10}}Q)$$
for all cubes $Q\subset \R^d$.
\end{lemma}

In a sense, Lemma \ref{fonam**}
acts as a substitute of the usual reverse H\"older
inequality for the classical $A_p$ weights. Its proof will follow by
a self improvement argument in the same spirit as the
proof of the reverse H\"older inequality for the $A_p$ weights.

Given $h\in\Z$ and $f\in L^1_{loc}(\mu)$, we denote
$$N^h f(x) = \sup_{k\geq h} S_k|f|(x).$$

The next technical result concentrates the main steps of the proof of
Lemma \ref{fonam**}.

\begin{lemma}\label{fonam'}
Let $S = \sup_Q \sigma(\frac{11}{10}Q)^{-1}\int |
N(\sigma\chi_Q)|^p\, w\,d\mu$, where the supremum is taken over all
cubes $Q\subset\R^d$. Assume that
$$\int |S_k(\sigma\,\chi_Q)|^p\, w\,d\mu \leq C\,\sigma(Q)$$
for all cubes $Q\subset \R^d$ uniformly on $k\in\Z$. Then,
for all $\ve>0$, there exists some constant $C_\ve$ such that for any
$\mu$-$\sigma$-$(2,\beta)$-doubling cube $Q\in\AD_{h}$,
\begin{equation}  \label{mmm1}
\int_Q |N^h(\sigma\,\chi_Q)|^p\, w\,d\mu \leq (C_\ve + \ve\,S)\,\sigma(Q).
\end{equation}
\end{lemma}

\begin{proof}
{\bf The construction.}
Let $Q_0$ be some fixed $\mu$-$\sigma$-$(2,\beta)$-doubling cube, with
$Q_0\in \AD_{h}$. We also denote
$\lambda_0:= m_{Q_0} \sigma$. We will show that \rf{mmm1} holds for $Q_0$.
To this end, by an inductive argument, for each $k\geq 1$ we will construct a
sequence of $\mu$-$\sigma$-doubling cubes $\{Q^k_i\}_i$.

First we will show how the cubes $\{Q^1_i\}$ are obtained. Let
\begin{equation} \label{mmm2}
\Omega_0 = \{N^{h+20}\sigma(x)>K\lambda_0\},
\end{equation}
where $K$ is some big positive constant which will be fixed below.
By Lemma \ref{open}, this set is open.
We consider some Whitney decomposition
$\Omega_0 = \bigcup_i R_i^1$, where $R_i^1$ are dyadic cubes
with disjoint interiors.

Let us check that $Q_0\setminus\Omega_0\neq
\varnothing$. If $Q_0\subset \Omega_0$, then for all $x\in Q_0\cap\supp(\mu)$
there exists some cube $Q_x$ centered at $x$, with
$Q_x\in\AD_{+\infty,h+19}$ with $m_{Q_x} \sigma>C\,K\,\lambda_0$ (where $C>0$
is some fixed constant). Since $Q_x\in\AD_{+\infty,h+19}$, we have
$\ell(Q_x) \leq\ell(Q_0)/10$. By Besicovitch's Covering Theorem, there
exists some covering $Q_0\subset \bigcup_i Q_{x_i}$ with finite
overlap. Using that $Q_0$ is $\sigma$-doubling, we obtain
\begin{eqnarray*}
\int_{Q_0} \sigma\,d\mu \geq C^{-1} \int_{2Q_0} \sigma\,d\mu & \geq & C^{-1}\sum_i
\int_{Q_{x_i}} \sigma\,d\mu \\
& \geq & C^{-1}\,K\,\lambda_0\sum_i \mu(Q_{x_i}) \geq
C^{-1}\,K\,\lambda_0\,\mu(Q_0).
\end{eqnarray*}
Therefore, $m_{Q_0}\sigma \geq C^{-1}\,K\,\lambda_0$, which is a contradiction if
$K$ is big enough.

Since $Q_0\setminus\Omega_0\neq\varnothing$, by the properties of the Whitney
covering, we have $\ell(R_j^1)\leq C_{20}\ell(Q_0)$ for any Whitney cube
$R^1_j$ such that $R^1_j\cap Q_0\neq \varnothing$. Moreover,
subdividing the Whitney cubes if necessary, we may assume that
$C_{20}\leq1/10$.

Let $g^1_j\in\Z$ be such that $R^1_j \in \AD_{g^1_j}$.
Observe that if $R^1_j\cap Q_0\neq\varnothing$, then
$R^1_j\subset\frac{3}{2}Q_0$, and so $g^1_j\geq h-2$.
For $x\in R^1_j\cap\supp(\mu)$, we consider some
$\mu$-$\sigma$-$(100,\beta)$-doubling cube
$Q_x^1\in\AD_{g^1_j+16}$, with $\beta$ given by Lemma \ref{dobsim}.
Now we take a Besicovitch's covering of $Q_0\cap\Omega_0$ with this type
of cubes: $Q_0\cap \Omega_0\subset\bigcup_{i\in I_1} Q^1_i$, and we define
$A_1:= \bigcup_{i\in I_1} Q^1_i$. Notice that, for each $i$, $10Q^1_i\subset
\frac{3}{2}Q_0$, because all the Whitney cubes intersecting $Q_0$ have
side length $\leq \ell(Q_0)/10$. In particular, we have
$A_1\subset \frac{3}{2}Q_0$. For each $i\in I_1$, let $h^1_i\in\Z$ be such
that $Q^1_i\in\AD_{h^1_i}$. If $Q^1_i$ is centered at some
point in $R^1_j$, then $h^1_i=g^1_j+16\geq h+14$.
This finishes the step $k=1$ of the construction.

Suppose now that the cubes $\{Q^k_i\}_{i\in I_k}$ (which are
$\mu$-$\sigma$-$(100,\beta)$-doubling, with $10Q^i_k\subset
\frac{3}{2}Q_0$, and have finite overlap) have already been
constructed. Let us see how the cubes $\{Q^{k+1}_i\}_{i\in I_{k+1}}$
are obtained. For each fixed cube $Q^k_i$ we repeat the arguments applied to
$Q_0$. We denote $\lambda^k_i=m_{Q^k_i}\sigma$ and let $h^k_i\in\Z$ be such
that $Q^k_i\in\AD_{h^k_i}$. We consider the open set
$\Omega^k_i = \{N^{h^k_i+20}\sigma(x)>K\lambda^k_i\}$, and a
decomposition of it into Whitney dyadic cubes with disjoint interiors:
$\Omega^k_i=\bigcup_j R^{k+1}_j$. Arguing as in the case of
$Q_0$, we  deduce $Q^k_i\setminus\Omega^k_i\neq
\varnothing$, and if $R^k_j\cap Q^k_i\neq \varnothing$, then
$\ell(R^k_j)\leq\ell(Q^k_i)/10$.
Given $g^{k+1}_j\in\Z$ such that $R^{k+1}_j \in
\AD_{g^{k+1}_j}$, for
$x\in R^{k+1}_j$, we consider  some $\mu$-$\sigma$-$(100,\beta)$-doubling
cube $Q_{i,x}^{k+1}\in \AD_{g^{k+1}_{j}+16}$.

It may happen that the union $\bigcup_{i\in I_k}
(\Omega^k_i\cap Q^k_i)$ is not
pairwise disjoint, and so for a fixed $x\in\bigcup_{i\in I_k}
(\Omega^k_i\cap Q^k_i)$
there are several indices $i$ such that $Q_{i,x}^{k+1}$ is defined.
In any case, for each  $x\in\bigcup_i (\Omega^k_i\cap Q^k_i)$ we
choose $Q_x^{k+1}:=Q_{i,x}^{k+1}$ with $i$ so that $x\in\Omega^k_i\cap
Q^k_i$ (no matter which $i$).
Now we take a Besicovitch covering of $\bigcup_i (\Omega^k_i\cap
Q^k_i)$ with cubes of the type $Q_x^{k+1}$. So we have
$\bigcup_{i\in I_k} (\Omega^k_i\cap Q^k_i) \subset \bigcup_{j\in
I_{k+1}} Q^{k+1}_j$,
and the cubes $Q^{k+1}_j$ have bounded overlap. Moreover, for each
$j\in I_{k+1}$ there exists some $i$ such that
$10Q^{k+1}_j\subset\frac{3}{2}Q^k_i\subset\frac{3}{2}Q_0$.
We define $A_{k+1}:= \bigcup_{j\in I_{k+1}} Q^{k+1}_j$, and we
denote by $h^{k+1}_j$ the integer such that $Q^{k+1}_j\in\AD_{
h^{k+1}_j}$.

% ****************************************************************************

\vspace{3mm}
{\bf The first step to estimate $\int_{Q_0}|N^{h}\sigma|^p\, w\,d\mu$.}
We want to show that given any $\ve>0$, if $K$ is big enough, then
\begin{equation} \label{mmm10}
\int_{Q_0} |N^{h}(\sigma\chi_{Q_0})|^p\, w\,d\mu
\leq (C_\ve + \ve\,S)\sum_{k=0}^\infty \sigma(A_k).
\end{equation}
We will prove this estimate inductively. First we deal with the
case $k=0$.
We have
\begin{eqnarray} \label{mmm10.1}
\int_{Q_0}\! |N^{h}(\sigma\chi_{Q_0}|^p\,w\,d\mu\! & \leq &\! \int_{Q_0}\!
\sum_{k=h}^{h+19} |S_k(\sigma\chi_{Q_0})|^p\,w\,d\mu +
\int_{Q_0}|N^{h+20}\sigma|^p\,w\,d\mu \nonumber \\
& \leq & C\,\sigma(Q_0) + \int_{Q_0}|N^{h+20}\sigma|^p\,w\,d\mu.
\end{eqnarray}

Given some small constant $\ve>0$, let $B_0=\{x\in
Q_0:\,S_{h+3}\sigma(x)\leq \ve\lambda_0\}$. Let us see that $\sigma(B_0)$ is
small. By Lemma \ref{guai'}, for all $x\in B_0$ there exists some
$\mu$-doubling cube $P_x\in\AD_{+\infty,h+2}$ centered at $x$ such that
$m_{2P_x}\sigma\leq C\ve\lambda_0$. We consider a Besicovitch's
covering of $B_0$ with this type of cubes. That is, $B_0\subset\bigcup_i
P_{x_i}$, with $\sum_i \chi_{P_{x_i}}\leq C$. We have
\begin{eqnarray*}
\sum_i \sigma(2P_{x_i}) & \leq & C\ve\lambda_0 \sum_i\mu(2P_{x_i}) \, \leq\,
C\ve\lambda_0 \sum_i \mu(P_{x_i}) \\
& \leq & C\ve\lambda_0\,\mu(2Q_0) \, \leq \, C\ve\lambda_0\,\mu(Q_0) \, =\,
C\ve\,\sigma(Q_0).
\end{eqnarray*}
In particular, we deduce $\sigma(B_0) \leq C\ve\,\sigma(Q_0)$.
Then we obtain
\begin{eqnarray} \label{mmm*4}
\int_{B_0} |N^{h+20}\sigma|^p\,w\,d\mu & \leq & \sum_i \int_{P_{x_i}}
|N^{h+20}\sigma|^p\,w\,d\mu \nonumber \\
& = & \sum_i \int_{P_{x_i}}
|N^{h+20}(\sigma\chi_{{\ts\frac{3}{2}}P_{x_i}})|^p\,w\,d\mu \nonumber\\
& \leq & S\sum_i\sigma({\ts\frac{11}{10}\frac{3}{2}}P_{x_i}) \,
\leq \, S\sum_i \sigma(2P_{x_i}) \, \leq\, C\ve S\,\sigma(Q_0).
\end{eqnarray}

Now we have to estimate $\int_{Q_0\setminus
  B_0}|N^{h+20}\sigma|^p\,w\,d\mu$. Given $x\in R^1_j\subset \Omega_0$,
let $x'\in\partial \Omega_0$ be such that
$|x-x'|=\dist(x,\R^d\setminus\Omega)$.
From Lemma \ref{perprmax}, we derive the following maximum principle:
\begin{equation} \label{mmm*1}
N^{h+25}(\sigma\chi_{\R^d\setminus 2R^1_j})(x) \leq C_{21}\,N^{h+20}\sigma(x')
\leq C_{21}K\lambda_0,
\end{equation}
where $C_{21}>1$ is some fixed constant depending on $C_0,n,d$.
Let us see that if
$N^{h+25}\sigma(x) > 2C_{21}K\lambda_0$, then
\begin{equation} \label{mmm*2}
N^{h+25}\sigma(x) \leq \max\Bigl( 2\max_{g^1_j-2\leq t
 \leq g^1_j+4} S_t(\sigma\chi_{2R^1_j})(x),\,
N^{g^1_j+5}(\sigma\chi_{2R^1_j})(x) \Bigr).
\end{equation}
Indeed, we have
$$N^{h+25}\sigma(x) \leq \max\Bigl(\max_{h+25\leq t
 \leq g^1_j+4} S_t\sigma(x),\, N^{g^1_j+5}\sigma(x) \Bigr),$$
(with equality if $h+25\leq g^1_j+5$).
If $N^{h+25}\sigma(x) \leq N^{g^1_j+5}\sigma(x)$, then \rf{mmm*2}
follows from the fact that
$N^{g^1_j+5}\sigma(x) =
N^{g^1_j+5}(\sigma\chi_{2R^1_j})(x)$.
If $N^{h+25}\sigma(x) =
S_{t_0}\sigma(x)$ for some $t_0$ with $h+25\leq t_0\leq g^1_j+4$, then
$S_{t_0}\sigma(x)>2C_{21}K\lambda_0$, and so
$$S_{t_0}(\sigma\chi_{2R^1_j})(x) \geq S_{t_0}\sigma(x) -
N^{h+25}(\sigma\chi_{\R^d\setminus2R^1_j})(x) \geq
\frac{1}{2}\,S_{t_0}\sigma(x),$$
by \rf{mmm*1}. Thus,
$$N^{h+25}\sigma(x) \leq 2S_{t_0}(\sigma\chi_{2R^1_j})(x) \leq 2
\max_{h+25\leq t \leq g^1_j+4} S_t(\sigma\chi_{2R^1_j})(x).$$
Moreover, it is easily checked that, for $t\leq g^1_j-2$ (and $x\in
R^1_j$), we have $S_t(\sigma\chi_{2R^1_j})(x)\leq
S_{g^1_j-2}(\sigma\chi_{2R^1_j})(x)$. Therefore, \rf{mmm*2} holds in
any case.

We denote $D_0:=\{x\in Q_0:\,N^{h+25}\sigma(x)
>2C_{21}K\lambda_0\}$. Notice that $D_0\subset\Omega_0\cap Q_0\subset
A_1$. We have
\begin{eqnarray*}
\int_{Q_0\setminus B_0}\!\!|N^{h+20}\sigma|^p\,w\,d\mu & \leq & \!
\sum_{t=h+20}^{h+24} \int_{Q_0\setminus B_0}\!\! |S_t\sigma|^p\,w \,d\mu +
\int_{Q_0\setminus B_0}\!\!|N^{h+25}\sigma|^p\,w\,d\mu  \\
& \leq & C\,\sigma(Q_0) +\int_{Q_0\setminus
  B_0}|N^{h+25}\sigma|^p\,w\,d\mu,
\end{eqnarray*}
where we have used that $S_t\sigma(x) = S_t(\sigma\chi_{2Q_0})(x)$ if
$t=20,\ldots,24$ and $x\in Q_0$.
Now we write
$$\int_{Q_0\setminus B_0}|N^{h+25}\sigma|^p\,w\,d\mu =
\int_{Q_0\setminus (B_0\cup D_0)} + \int_{D_0\setminus B_0} =: I +
II.$$

First we will estimate $I$. For $x\in Q_0\setminus(B_0\cup D_0)$, we
have
$$N^{h+25}\sigma(x) \leq CK\lambda_0 \leq
CK\ve^{-1}\,S_{h+3}\sigma(x).$$
Therefore,
\begin{eqnarray*}
I & = & \int_{Q_0\setminus (B_0\cup D_0)}|N^{h+25}\sigma|^p\,w\,d\mu
\, \leq \,CK^p\ve^{-p} \int_{Q_0} |S_{h+3}\sigma|^p\,w\,d\mu \\
& \leq & CK^p\ve^{-p}\sigma(2Q_0) \, \leq \, CK^p\ve^{-p}\sigma(Q_0),
\end{eqnarray*}
where we have used that $S_{h+3}\sigma(x) =
S_{h+3}(\sigma\chi_{2Q_0})(x)$.

It remains to estimate $II$. Given $x\in R^1_j\cap(D_0\setminus B_0)$, by
\rf{mmm*2} we get
% $$N^{h+25}\sigma(x) \leq \max\Bigl( 2\max_{g^1_j-2\leq t
% \leq g^1_j+39} S_t(\sigma\chi_{2R^1_j})(x),\,
% N^{g^1_j+40}(\sigma\chi_{2R^1_j})(x) \Bigr).$$
% Then we get
\begin{eqnarray*}
\int_{R^1_j\cap(D_0\setminus B_0)} |N^{h+25}\sigma|^p\,w\,d\mu
& \leq & C\sum_{t=g^1_j-2}^{g^1_j+39}
\int |S_t(\sigma\chi_{2R^1_j})|^p\,w\,d\mu \\
&& \mbox{} +
\int_{R^1_j\cap(D_0\setminus B_0)} |N^{g^1_j+40}\sigma|^p\,w\,d\mu \\
& \leq & C\,\sigma(2R^1_j) +
\int_{R^1_j\cap(D_0\setminus B_0)} |N^{g^1_j+40}\sigma|^p\,w\, d\mu.
\end{eqnarray*}
Given $k\geq1$, for $x\in A_k$, we denote
$H^k_x := \max\{h^k_i:\,i\in I_k,\,x\in Q^k_i\}$.
It is easily seen that if $x\in R^1_j\cap A_1$, then
$H^1_x +20 \leq g^1_j+40$.
Then, summing over all the cubes $R^1_j\subset \Omega_0$ such that $R^1_j\cap
Q_0\neq\varnothing$, due to the finite overlap of the cubes $2R^1_j$, we
obtain
\begin{equation}\label{nn*21}
\int_{D_0\setminus B_0} |N^{h+25}\sigma|^p\,w\,d\mu
\leq C\,\sigma(2Q_0) + \int_{A_1} |N^{H^1_x+20}\sigma(x)|^p\,w(x)\,d\mu(x).
\end{equation}
So we have shown that
\begin{equation} \label{mmm*10}
  \int_{Q_0}|N^h(\sigma\chi_{Q_0})|^p\,w\,d\mu \leq
(C_{22} + C_{23}\ve S)\,\sigma(Q_0) + \int_{A_1}
|N^{H^1_x+20}\sigma|^p\,w\,d\mu,
\end{equation}
with $C_{22}$, but not $C_{23}$, depending on $K$ and $\ve$.

% ***************************************************************************

\vspace{3mm}
{\bf The $k$-th step to estimate  $\int_{Q_0}|N^{h}\sigma|^p\, w\,d\mu$.}
Now we will show that for any $k\geq1$,
\begin{equation} \label{mmm*11}
\int_{A_k}  |N^{H^k_x+20}\sigma|^p\,w\,d\mu \leq (C_{22}'+
C_{23}'\ve S)\,\sigma(A_k) + \int_{A_{k+1}}
|N^{H^{k+1}_x+20}\sigma|^p\,w\,d\mu,
\end{equation}
with $C_{22}'$, but not $C_{23}'$, depending on $K$ and $\ve$.
The arguments to prove \rf{mmm*11} are similar to the ones we have used
to obtain \rf{mmm*10}, although a little more involved because the
cubes $\{Q^k_i\}_{i\in I_k}$ are non pairwise disjoint.

For each $i\in I_k$ we define $B^k_i = \{x\in Q^k_i:\,
S_{h^k_i+3}\sigma(x) \leq \ve\lambda^k_i\}$. Arguing as in \rf{mmm*4},
we deduce
$$\int_{B^k_i}|N^{h^k_i+20}\sigma|^p\,w\,d\mu \leq
C\ve S\,\sigma(Q^k_i).$$
We denote $B_k=\bigcup_{i\in I_k} B^k_i$. Using the definition of
$H^k_x$, we obtain
\begin{eqnarray} \label{ames}
\int_{B_k} |N^{H^k_x+20}\sigma(x)|^p\,w(x)\,d\mu(x) & \leq & \sum_i
\int_{B^k_i} |N^{H^k_x+20}\sigma(x)|^p\,w(x)\,d\mu(x) \nonumber \\
& \leq & \int_{B^k_i} |N^{h^k_i+20}\sigma|^p\,w\,d\mu \nonumber\\
& \leq & C\ve S\sum_i \sigma(Q^k_i) \,\leq \, C\ve S\sigma(A_k).
\end{eqnarray}

To estimate $\int_{A_k\setminus B_k} |N^{H^k_x+20}\sigma|^p\,w\,d\mu$,
we need to introduce some additional notation.
Assume $I_k=\{1,2,3,\ldots\}$. We denote $I_{k,t}:=\{i\in I_k:Q^k_i\in
\AD_{t}\}$. We set
$$\wh{Q}^k_i :=Q^k_i\setminus \Bigl( \bigcup_{l \in I_{k,t},\,t>h^k_i}
Q^k_l \,\cup\!  \bigcup_{l \in I_{k,h^k_i},\,l<i} Q^k_l \Bigr).$$
It is easily checked that the sets $\wh{Q}^k_i$, $i\in I_k$,  are
pairwise disjoint, that $\bigcup_{i\in I_k} \wh{Q}^k_i = \bigcup_{i\in
I_k} Q^k_i = A_k$, and moreover that if $x\in\wh{Q}^k_i$, then
$H^k_x=h^k_i$. We have
\begin{eqnarray} \label{mb1}
\lefteqn{\int_{A_k\setminus B_k}|N^{H^k_x+20}\sigma|^p\,w\,d\mu \, = \,
\sum_{i\in I_k}
\int_{\wh{Q}^k_i\setminus B_k}|N^{h^k_i+20}\sigma|^p\,w\,d\mu} &&
\nonumber\\
& \leq & \sum_{i\in I_k}
\sum_{t=h^k_i+20}^{h^k_i+24} \int_{\wh{Q}^k_i\setminus B_k}
|S_t\sigma|^p\,w\,d\mu + \sum_{i\in I_k} \int_{\wh{Q}^k_i\setminus B_k}
|N^{h^k_i+25}\sigma|^p\,w\,d\mu\nonumber\\
& \leq & C\sum_{i\in I_k}\sigma(2Q^k_i) + \int_{A_k\setminus B_k}
|N^{H^k_x+25}\sigma|^p\,w\,d\mu \nonumber\\
& \leq & C\sigma(A_k) + \int_{A_k\setminus B_k}
|N^{H^k_x+25}\sigma|^p\,w\,d\mu.
\end{eqnarray}

Now we set $D^k_i =\{x\in Q^k_i:\,N^{h^k_i+25}\sigma >
2C_{21}K\lambda^k_i\}$, and $D_k = \bigcup_{i\in I_k} D^k_i$.
For $x\in Q^k_i\setminus(D^k_i\cup B_k)$, we have
$$N^{H^k_x+25}\sigma(x) \leq N^{h^k_i+25}\sigma(x) \leq C
K\ve^{-1}\,S_{h^k_i+3}\sigma(x).$$
Therefore, operating as in the case $k=0$, we get
$$\int_{Q^k_i\setminus D^k_i\cup B_k}
|N^{H^k_x+25}\sigma(x)|^p\,w(x)\,d\mu(x) \leq CK^p\ve^{-p}\,\sigma(Q^k_i).$$
Summing over $i\in I_k$, we obtain
\begin{equation} \label{mb2}
\int_{A_k\setminus(B_k\cup D_k)}|N^{H^k_x+25}\sigma(x)|^p\,w(x)\,d\mu(x)
\leq  CK^p\ve^{-p}\,\sigma(A_k).
\end{equation}

Finally we deal with
$\int_{D_k\setminus B_k}|N^{H^k_x+25}\sigma(x)|^p\,w(x)\,d\mu(x)$. For
a fixed $k$,
let $\{R^{k+1}_j\}_{j\in J_{k+1}}$ be the collection of {\em all} the
Whitney cubes (originated from {\em all} the sets $\Omega^k_i$ $i\in I_k$)
such that if $R^{k+1}_j$ comes from $\Omega^k_i$, then $R^{k+1}_j\cap
Q^k_i\neq \varnothing$.
Assume $J_{k+1}=\{1,2,3,\ldots\}$. We denote $J_{k+1,t}:=\{j\in
J_{k+1}:R^{k+1}_j\in\AD_{t}\}$. We set
$$\wh{R}^{k+1}_j := {\ts\frac{3}{2}}R^{k+1}_j\setminus
\Bigl( \bigcup_{l \in J_{k+1,t},\,t>g^{k+1}_j} {\ts\frac{3}{2}}
R^{k+1}_l \,\cup\!  \bigcup_{l \in J_{k+1,g^{k+1}_j},\,l<j}
{\ts\frac{3}{2}}R^{k+1}_l \Bigr).$$
The sets $\wh{R}^{k+1}_j$, $j\in J_{k+1}$,  are
pairwise disjoint and
$$\bigcup_{j\in J_{k+1}} \wh{R}^{k+1}_j = \bigcup_{j\in
J_{k+1}} \frac{3}{2}R^{k+1}_j\supset A_{k+1}.$$
Moreover, it easily seen that if
$x\in\wh{R}^{k+1}_j$, then $g^{k+1}_j+40\geq H^{k+1}_x+20$, and so
$N^{g^{k+1}_j+40}\sigma(x) \leq N^{H^{k+1}_x+20}\sigma(x).$
If $\wh{R}^{k+1}_j$ is originated by $\Omega^k_i$, arguing as in the
case $k=0$, we deduce
\begin{eqnarray*}
\lefteqn{N^{H^k_x+25}\sigma(x) \, \leq \, N^{h^k_i+25}\sigma(x)}&&\\
& \leq & \max \Bigl(2\max_{g^{k+1}_{j}-2\leq t
\leq g^{k+1}_j+39} S_t(\sigma\chi_{2R^{k+1}_j})(x),\,
N^{g^{k+1}_j+40}(\sigma\chi_{2R^{k+1}_j})(x) \Bigr).
\end{eqnarray*}
Therefore,
\begin{eqnarray} \label{mb3}
\lefteqn{
\int_{D_k\setminus B_k} |N^{H^k_x+25}\sigma(x)|^p\,w(x)\,d\mu(x) \, = \,
\sum_{j\in J_{k+1}} \int_{\wh{R}^{k+1}_j\cap (D_k\setminus B_k)}}&& \nonumber\\
& \leq & \sum_{j\in J_{k+1}} \sum_{t=g^{k+1}_j-2}^{g^{k+1}_j+39} \int
|S_t(\sigma\chi_{2R^{k+1}_j})|^p\,w\,d\mu \nonumber\\
&& \mbox{}+ \sum_{j\in J_{k+1}}  \int_{\wh{R}^{k+1}_j\cap (D_k\setminus B_k)}
|N^{H^{k+1}_x+20}\sigma(x)|^p\,w(x)\,d\mu(x) \nonumber \\
& \leq & C\sum_{j\in J_{k+1}}\sigma(2R^{k+1}_j) + \int_{A_{k+1}}
|N^{H^{k+1}_x+20}\sigma(x)|^p\,w(x)\,d\mu(x).
\end{eqnarray}
In the following estimates the notation $R^{k+1}_j\sim Q^k_i$ means that
$R^{k+1}_j$ is a Whitney cube of $\Omega^k_i$:
\begin{eqnarray} \label{mb4}
\sum_{j\in J_{k+1}}
\sigma(2R^{k+1}_j) & = & \sum_{i\in I_k}\,\, \sum_{j\in J_{k+1}:
\,R^{k+1}_j\sim Q^k_i} \sigma(2R^{k+1}_j) \nonumber \\
& \leq & C\sum_{i\in I_k} \sigma(\Omega^k_i\cap 2Q^k_i)\, \leq \,
C\sum_{i\in I_k}\sigma(Q^k_i) \,\leq \,C\,\sigma(A_k).
\end{eqnarray}
By \rf{ames}, \rf{mb1}, \rf{mb2}, \rf{mb3} and \rf{mb4}, \rf{mmm*11} follows.

From \rf{mmm*10} and \rf{mmm*11}, we get
\begin{eqnarray} \label{mb10}
\int_{Q_0}N(\sigma\chi_{Q_0})|^p\,w\,d\mu & \leq & (C + C\ve S)
  \sum_{k=0}^\infty \sigma(A_k) \nonumber\\
&& \mbox{} + \limsup_{k\to\infty}\int_{A_k}
|N^{H^{k}_x+20}\sigma(x)|^p\,w(x)\,d\mu(x).
\end{eqnarray}
This is the same as \rf{mmm10}, except for the last term on right hand
side. However, we will see below that this term equals $0$.

% ****************************************************************************

\vspace{3mm}
{\bf The estimate of $\sum_k \sigma(A_k)$.} We are going to prove that
\begin{equation} \label{mmm50}
\sum_{k=0}^\infty \sigma(A_k) \leq C\,\sigma(Q_0).
\end{equation}
We denote $\wt{A}_k=\bigcup_{i\in I_k} 2Q^k_i$. It is easily seen that
$\wt{A}_{k+1} \subset \wt{A}_k$ for all $k$ (this is the main
advantage of $\wt{A}_k$ over $A_k$). We will show that there exists
some positive constant $\tau_0<1$ such that
\begin{equation} \label{mmm51}
\sigma(\wt{A}_{k+2}) \leq \tau_0\,\sigma(\wt{A}_k)
\end{equation}
for all $k$. This implies \rf{mmm50}, because
$\wt{A}_0,\wt{A}_1\subset 2Q_0$ and $Q_0$ is $\sigma$-doubling.

For a fixed $k\geq1$, by the covering Lemma \ref{vitali}, there exists
some subfamily $\{2Q^k_j\}_{j\in I_k^0}\subset \{2Q^k_i\}_{i\in I_k}$
such that
\begin{itemize}
\item[(1)] $\wt{A}_k\subset \bigcup_{j\in I_k^0} 40Q^k_j$.
\item[(2)] $4Q^k_j\cap4Q^k_l = \varnothing$ if $j,l\in I_k^0$.
\item[(3)] If $j\in I_k^0$, $l\not\in I_k^0$, and
  $4Q^k_j\cap4Q^k_l\neq\varnothing$, then $\ell(Q^k_l) \leq
  10\ell(Q^k_j)$.
\end{itemize}
First, we will see that
\begin{equation} \label{mmm52}
\sigma(2Q^k_j\cap \wt{A}_{k+2}) \leq \tau_1\,\sigma(2Q^k_j) \qquad
\mbox{if $j\in I_k^0$,}
\end{equation}
for some fixed constant $0<\tau_1<1$. By Lemma \ref{cubspetits}
it is enough to show that, for each $x\in \frac{3}{2}Q^k_j\cap\supp(\mu)$,
there exists some cube
$P\in\AD_{h^k_j+4}$ centered at $x$ such that
$\mu(\wt{A}_{k+2}\cap P)\leq \delta_0\,\mu(P)$, with $\delta_0$
sufficiently small.

Let $2Q_s^{k+1}$ some cube which forms $\wt{A}_{k+1}$ such
that $2Q^{k+1}_s\cap 2Q^k_j\neq\varnothing$.
By our construction, there exists some cube $Q^k_t$ such that
$10Q^{k+1}_s\subset\frac{3}{2}Q^k_t$, so that $Q^{k+1}_s$ comes from
$\Omega^k_t$. Because of the property (3) of the
covering, we have $\ell(Q^k_t)\leq10\ell(Q^k_j)$. Therefore, $2Q^k_t\in
\AD_{+\infty,h^k_j-3}$, which implies $Q^{k+1}_s \in \AD_{+\infty,h^k_j+7}$ and
$2Q^{k+1}_s \in \AD_{+\infty,h^k_j+6}$.

Let $P\in\AD_{h^k_j+4}$ be some $\mu$-doubling cube whose center is
in $\frac{3}{2}Q^k_j$ (which implies $P\subset 2Q^k_j$).
Let $S_P$ be the set of indices $s$ such that
$2Q^{k+1}_s\cap P\neq\varnothing$. We have $\ell(Q^{k+1}_s)
\ll\ell(P)$ for $s\in S_P$, because
$2Q^{k+1}_s\in\AD_{+\infty,h^k_j+6}$ (since $2Q^{k+1}_s \cap
2Q^k_j\neq \varnothing$). Thus, $2Q^{k+1}_s \subset 2P$.
From our construction, we deduce
$$\mu(\wt{A}_{k+2}\cap P) \leq \mu\Bigl(\bigcup_{s\in S_P} (\Omega^{k+1}_s
\cap 2Q^{k+1}_s)\Bigr) \leq \sum_{s\in S_P} \mu(\Omega^{k+1}_s
\cap 2Q^{k+1}_s).$$
Since
$N^{h^{k+1}_s+20}\sigma(x) = N^{h^{k+1}_s+20}(\sigma\chi_{3Q^{k+1}_s})(x)$
for $x\in Q^{k+1}_s$, by the weak $(1,1)$ boundedness of
$N^{h^{k+1}_s+20}$, and by the $\sigma$-doubling property of
$Q^{k+1}_s$, we obtain
$$\mu(\Omega^{k+1}_s\cap 2Q^{k+1}_s) \leq
\frac{C\,\sigma(3Q^{k+1}_s)}{K\,\lambda^{k+1}_s} \leq
\frac{C}{K}\,\mu(Q^{k+1}_s).$$
Thus, by the finite overlap of the cubes $Q^{k+1}_s$ and the fact that
$P$ is $\mu$-doubling,
$$\mu(\wt{A}_{k+2}\cap P) \leq \frac{C}{K}\sum_{s\in S_P}
\mu(Q^{k+1}_s) \leq \frac{C}{K} \,\mu(2P) \leq\frac{C}{K} \,\mu(P)
=: \delta_0\mu(P).$$
Since we may choose $K$ as big as we want, $\delta_0$ can be taken
arbitrarily small, and \rf{mmm52} holds.

Let us see that \rf{mmm51} follows from \rf{mmm52}. We denote
$\wt{A}_{k,0} = \bigcup_{j\in I_k^0} 2Q^k_j$.
Since the cubes $2Q^k_j$, $j\in I_k^0$, are disjoint, \rf{mmm52}
implies
$\sigma(\wt{A}_{k,0}\cap \wt{A}_{k+2}) \leq \tau_1\,\sigma(\wt{A}_{k,0})$.
By the property (1) of the covering and the fact that $Q^k_j$ is
$(100,\beta)$-$\sigma$-doubling, we have
$$\sigma(\wt{A}_{k,0}) = \sum_{j\in I_k^0} \sigma(2Q^k_j) \geq
C_{24}^{-1} \sum_{j\in I_k^0} \sigma(40Q^k_j) \geq C_{24}^{-1}
\sigma(\wt{A}_k).$$
Then,
$$\sigma(\wt{A}_k\setminus \wt{A}_{k+2}) \geq
\sigma(\wt{A}_{k,0}\setminus \wt{A}_{k+2})
\geq (1-\tau_1)\,\sigma(\wt{A}_{k,0}) \geq
(1-\tau_1)\,C_{24}^{-1}\,\sigma(\wt{A}_k).$$
Therefore,
$$\sigma(\wt{A}_k\cap \wt{A}_{k+2}) \leq (1-(1-\tau_1)\,C_{24}^{-1})\,
\sigma(\wt{A}_k) =: \tau_0\,\sigma(\wt{A}_k).$$

% ****************************************************************************

\vspace{3mm}
{\bf The end of the proof.} We only need to prove the lemma for
$S<+\infty$. For each $k\geq1$ we have
\begin{eqnarray*}
\int_{A_k} |N^{H^{k}_x+20}\sigma(x)|^p\,w(x)\,d\mu(x) & \leq &
\sum_{i\in I_k} \int_{Q^k_i}|N^{h^k_i+20}\sigma(x)|^p\,w(x)\,d\mu(x)\\
& \leq & CS\sum_i \sigma(Q^k_i) \,\leq\, CS\, \sigma(A_k).
\end{eqnarray*}
From \rf{mmm50} we deduce that $\sigma(A_k)\to0$ as $k\to\infty$, and
then the integral on the left
hand side above tends to $0$ as $k\to\infty$. Now the lemma
follows from \rf{mb10} and \rf{mmm50}.
\end{proof}

% ****************************************************************************
% ****************************************************************************

\begin{proof}[{\bf Proof of Lemma \ref{fonam**}}]
Let $Q$ be some cube with $Q\in \AD_{h}$ and $x_0\in Q\cap\supp(\mu)$. We
write
$$\int |N(\sigma\chi_Q)|^p\,w\,d\mu = \int_{\frac{21}{20} Q} +
\int_{Q_{x_0,h-4} \setminus\frac{21}{20} Q} + \int_{\R^d\setminus
Q_{x_0,h-4}} =: I + II + III.$$
First we will estimate the integral $I$. For each $x\in
\frac{21}{20} Q\cap
\supp(\mu)$, let $P_{x}$ be some $\mu$-$\sigma$-$(4,\beta)$-doubling cube with
$P_x\in\AD_{h+10}$.
Notice that for each $y\in P_{x}$ and $k\geq h+15$, we have
$\supp(s_k(y,\cdot))\subset 2P_{x}$. Thus by Lemma \ref{fonam'}, if we
denote $C_S := C_\ve + \ve S$, we get
$$\int_{P_x} |N^{h+15}\sigma|^p\,w\,d\mu =
\int_{P_x} |N^{h+15}(\sigma\chi_{2P_x})|^p\,w\,d\mu \leq C_S\,\sigma(2P_x)
\leq C\,C_S\,\sigma(P_x).$$
By Besicovitch's Covering Theorem,
there exists some subfamily of cubes $\{P_{x_i}\}_i\subset \{P_x\}_x$
which covers $\frac{21}{20} Q\cap\supp(\mu)$ with finite overlap.
Since $\ell(P_{x_i})\ll\ell(Q)$, we have
$P_{x_i}\subset \frac{11}{10}Q$. Then we obtain
\begin{eqnarray*}
\int_{\frac{21}{20}Q} |N^{h+15}\sigma|^p\,w\,d\mu &
\leq & \sum_i \int_{P_{x_i}} |N^{h+15}\sigma|^p\,w\,d\mu \\
& \leq & C\,C_S\sum_i \sigma(P_{x_i})
\,\leq\, C\,C_S\,\sigma(\ts\frac{11}{10}Q).
\end{eqnarray*}
It is easily seen that, for all $y\in \frac{21}{20}Q$,
$N(\sigma\chi_Q)(y) \leq C\,N^{h-2}(\sigma\chi_Q)(y)$.
Therefore,
\begin{eqnarray*}
I & \leq & C\int_{\frac{21}{20}Q} |N^{h-2}(\sigma\chi_Q)|^p\,w\,d\mu \\
& \leq & C\int_{\frac{21}{20}Q} \bigl|\sum_{k=h-2}^{h+14} S_k(
\chi_Q\sigma)\bigr|^p\,w\,d\mu + \int_{\frac{21}{20}Q}
|N^{h+15}\sigma|^p\,w\,d\mu \\
& \leq & C(1+C_S)\sigma(\ts\frac{11}{10}Q).
\end{eqnarray*}

Now we turn our attention to the integral $II$. For $y\in Q_{x_0,h-4}
\setminus \frac{21}{20}Q$,
$$N(\sigma\chi_Q)(y) \leq \frac{C\,\sigma(Q)}{|y-x_0|^n}
\leq C \sum_{k=h-6}^{h+3} S_k(\sigma\chi_Q)(y).$$
Thus,
$$II \leq \int |\sum_{k=h-6}^{h+3} S_k(\sigma\chi_Q)|^p\,w\,d\mu \leq
C\,\sigma(Q).$$

Finally we deal with $III$. For $k\leq h-4$ and $y\in
Q_{x_0,k-1}\setminus Q_{x_0,k}$, we have
$$N(\sigma\chi_Q)(y) \leq C\,\frac{\sigma(Q)}{|y-x_0|^n} \leq C\,
\frac{\sigma(Q)}{\sigma(Q_{x_0,k+1})}\, \sum_{j=k-3}^{k+2}
S_{j} (\sigma\chi_{Q_{x_0,k+1}})(y).$$
Thus,
\begin{eqnarray*}
\int_{Q_{x_0,k}\setminus Q_{x_0,k-1}}\!\!\!\! |N(\sigma\chi_Q)|^p\,w\,d\mu
& \leq &  \frac{C\sigma(Q)^p}{\sigma(Q_{x_0,k+1})^p}
\int \Bigl|\sum_{j=k-3}^{k+2}\! S_{j} (\sigma
\chi_{Q_{x_0,k+1}})\Bigr|^p w\,d\mu \\
& \leq & \frac{C\sigma(Q)^p}{\sigma(Q_{x_0,k+1})^{p-1}}.
\end{eqnarray*}
From Lemma \ref{sigmenor}, we deduce
$\sigma(Q) \leq \sigma(Q_{x_0,h-1})
\leq \eta^{h-k-2} \sigma(Q_{x_0,k+1})$
for some fixed constant $\eta$ with $0<\eta<1$.
Therefore,
\begin{eqnarray*}
III & = & \sum_{k=-\infty}^{h-4}
\int_{Q_{x_0,k-1}\setminus Q_{x_0,k}} |N(\sigma\chi_Q)|^p\,w\,d\mu \\
%& \leq & C\,\sigma(Q) \sum_{k=-\infty}^{h-4}
% \frac{\sigma(Q)^{p-1}}{\sigma(Q_{x_0,k+1})^{p-1}}\\
& \leq & C\,\sigma(Q) \sum_{k=-\infty}^{h-4} \eta^{(p-1)(h-k-2)}
\, \leq \,C\,\sigma(Q).
\end{eqnarray*}

So we have
$$\int |N(\sigma\chi_Q)|^p\,w\,d\mu \leq
C(1+C_S)\,\sigma({\ts\frac{11}{10}Q})
= C_{25}\,(1+C_\ve+\ve S)\,\sigma({\ts\frac{11}{10}Q}).$$
Choosing $\ve\leq 1/(2C_{25})$ and taking the supremum over all the
cubes $Q$, we get $S\leq C_{25}(1+C_\ve) + \frac{1}{2}S$. Thus
$S\leq 2C_{25}(1+C_\ve)$ if $S<+\infty$.

One way to avoid the assumption $S<+\infty$ would be
to work with ``truncated'' operators of type
$N^{h,l}f:=\sup_{h\leq k\leq l}S_k|f|$
in Lemma \ref{fonam'}, instead of $N^h$; and also to consider
a truncated version of $S$ in \rf{mmm1}, etc. The technical details are
left for the reader.
\end{proof}

% ****************************************************************************
% ****************************************************************************
% ****************************************************************************

\section{Boundedness of $N$ on $L^p(w)$}

The implication (e) $\Rightarrow$ (c) of Lemma \ref{main} follows
from Lemma \ref{fonam**} and the following result.

\begin{lemma}  \label{dif2}
If for any $k\in\Z$ and any cube $Q$,
\begin{equation} \label{cond1}
\int N(\sigma\,\chi_Q)^p\, w\,d\mu \leq C\,\sigma({\ts \frac{11}{10}}Q)
\end{equation}
and
\begin{equation} \label{cond2}
\int S_k(w\,\chi_Q)^{p'}\, \sigma\,d\mu \leq C\,w(Q),
\end{equation}
with $C$ independent of $k$ and $Q$, then $N$ is bounded on $L^p(w)$.
\end{lemma}

The proof of this lemma is inspired by the techniques used by Sawyer
\cite{Sawyer} to obtain two weight norm inequalities for fractional
integrals. In our case, we have to overcome new difficulties which are
mainly due to the fact that the operator $N$ is not linear and it is
very far from behaving as a self adjoint operator, because it is a {\em
centered} maximal operator.

\begin{proof}
We will show that for some $\beta\geq0$ the operator
$T:= N + \beta\,M_R$ is bounded on $L^p(w)$.
Without loss of generality, we take $f\in L^1(\mu)$ non negative with
compact support.
Given some constant $\alpha>1$ close to $1$, for each $k\in\Z$, we denote
$$\Omega_k = \{x:\, Tf(x)>\alpha^k\}.$$
The precise value of $\alpha$ and $\beta$ will
be fixed below.
As in Lemma \ref{prmax}, we consider the Whitney decomposition
$\Omega_k = \bigcup_i Q^k_i$, where $Q^k_i$ are dyadic cubes with disjoint
interiors (the {\it Whitney cubes}).

Take some cube $Q^k_i\subset \Omega_k$ and $x\in Q^k_i \cap \Omega_{k+2}$.
Suppose that $m$ and $\beta$ are chosen in Lemma \ref{prmax} so that the
maximum principle \rf{eqprmax} holds with $\ve=\alpha-1$.
Then, we have
\begin{equation} \label{vv****}
T(f\,\chi_{\R^d\setminus U_m(Q^k_i)})(x) \leq (1+\ve)\,\alpha^k =
\alpha^{k+1},
\end{equation}
and so
\begin{equation}  \label{andso}
T(f\,\chi_{U_m(Q^k_i)})(x) \geq \alpha^{k+2} - \alpha^{k+1}=
\frac{\alpha-1}{\alpha}\,\alpha^{k+2}.
\end{equation}
Let $h\in \Z$ be such that $Q_i^k\in \AD_{h}$.
If for all $j$ with $h-M\leq j \leq h+M$ (where $M$ is some positive
big integer which will be chosen later) we have
$$S_j (f\chi_{U_m(Q^k_i)}) (x) \leq \delta\,\alpha^k,$$
where $\delta>0$ is another constant which we will choose below, then
we write $x\in B_k$ (i.e. $x$ is a ``bad point'').
Notice that $B_k \subset\Omega_{k+2}\subset \Omega_k$.

We will see that the set of bad points is quite small.
Indeed, we will prove that
\begin{equation}  \label{fonam}
w\biggl( \bigcup_{j\geq k} B_j \biggr) \leq \eta\, w(\Omega_k),
\end{equation}
where $0<\eta<1$ is some constant which depends on $\tau$ (from
the $Z_\infty$ property), $n$, $d$, but not on
$\beta$, $m$, $\alpha$, $M$.
We defer the proof of \rf{fonam}, which is one of the key points of our
argument, until Lemma \ref{lemaclau} below.

Let us denote $A_k=\bigcup_{j\geq k} B_j$. Now we have
\begin{eqnarray}  \label{sim2}
\int |Tf|^p\,w\,d\mu & = & \int_0^\infty p\,\lambda^{p-1}\,w(\Omega_\lambda)\,
d\lambda \nonumber\\
& \leq & \sum_{k\in\Z} p\,(\alpha^{k+1} - \alpha^k)\,\alpha^{(k+1)(p-1)}\,
w(\Omega_k) \nonumber\\
& = & p\,\alpha^{p-1}\,(\alpha-1) \sum_{k\in\Z} \alpha^{kp}\,
\bigl[w(\Omega_k\setminus A_{k-2}) + w(\Omega_k\cap A_{k-2})\bigr].
\end{eqnarray}
From \rf{fonam} we get
\begin{equation}  \label{eq*}
p\,\alpha^{p-1}\,(\alpha-1) \sum_{k\in\Z}
\alpha^{kp}\,w(\Omega_k\cap A_{k-2}) \leq
\eta \,p\,\alpha^{p-1}\,(\alpha-1) \sum_{k\in\Z} \alpha^{kp}\,w(\Omega_{k-2}).
\end{equation}
From calculations similar to the ones in \rf{sim2}, it follows
$$\int |Tf|^p\,w\,d\mu \geq p\,(\alpha-1)\,\alpha^{-3p} \sum_{k\in\Z}
\alpha^{kp}\,w(\Omega_{k-2}).$$
If we take $\alpha$ such that $\eta^{1/2}\alpha^{4p-1} = 1$, then the right
hand side of \rf{eq*} is bounded above by $\eta^{1/2}\int|Tf|^p\,w\,d\mu$,
and so
$$\int|Tf|^p\,w\,d\mu \leq (1-\eta^{1/2})^{-1}
\sum_{k\in\Z} \alpha^{kp}\,w(\Omega_k\setminus A_{k-2}).$$
Summing by parts we get
\begin{eqnarray*}
\int|Tf|^p\,w\,d\mu & \leq & C\,\sum_{k\in\Z} \alpha^{kp}\,
w(\Omega_{k+2}\setminus A_{k}) \\
& = & C\,\sum_{k\in\Z} \alpha^{kp}\,\bigl[w(\Omega_{k+2}\setminus A_k)
- w(\Omega_{k+3}\setminus A_{k+1})\bigr].
\end{eqnarray*}
Observe that if we assume $\int|Tf|^p\,w\,d\mu<\infty$, then
$$\sum_{k\in\Z} \alpha^{kp}\,w(\Omega_{k+2}\setminus A_k) <\infty\quad
\mbox{and} \quad\sum_{k\in\Z} \alpha^{kp}\,w(\Omega_{k+3}\setminus A_{k+1})
<\infty,$$
which implies that our summation by parts is right.
Since $A_{k+1}\subset A_k$, we have
$w(\Omega_{k+3}\setminus A_{k+1}) \geq w(\Omega_{k+3}\setminus A_{k}).$
Thus,
\begin{equation}  \label{eq**}
\int|Tf|^p\,w\,d\mu \leq C\,\sum_{k\in\Z} \alpha^{kp}\,
w((\Omega_{k+2}\setminus \Omega_{k+3})\setminus A_{k}).
\end{equation}

We denote $E^k_i = Q^k_i \cap
(\Omega_{k+2}\setminus \Omega_{k+3})\setminus A_{k}$ for all $(k,i)$.
To simplify notation, we also set $U^k_i = U_m(Q^k_i)$.
Given $h\in\Z$ such that $Q^k_i\in \AD_{h}$, we consider the operator
$$\sk = S_{h-M} + S_{h-M+1} + \cdots + S_{h+M}.$$
Since $E^k_i\subset \R^d\setminus A_k$, we obtain
\begin{eqnarray*}
w(E^k_i) & \leq & \delta^{-1}\alpha^{-k} \int_{E^k_i} \sk(\chi_{U^k_i}\,f)
\,w\,d\mu \\
& = & \delta^{-1}\alpha^{-k} \int_{U^k_i} f\,\skd(\chi_{E^k_i}\,w)\,d\mu \\
& = & \delta^{-1}\alpha^{-k} \left(\int_{U^k_i\setminus\Omega_{k+3}}
+ \int_{U^k_i\cap\Omega_{k+3}} \right) \, =\,
\delta^{-1}\alpha^{-k} (\sigma^k_i + \tau^k_i).
\end{eqnarray*}
From \rf{eq**} we get
\begin{eqnarray}  \label{eq123}
\int|Tf|^p\,w\,d\mu & \leq & C\,\sum_{k,i} \alpha^{kp}\,w(E^k_i) \nonumber \\
& = & C\, \left( \sum_{(k,i)\in E} + \sum_{(k,i)\in F} +
\sum_{(k,i)\in G}\right)  \cdot \alpha^{kp}\,w(E^k_i) \nonumber \\
& = & C\,(I + II + III),
\end{eqnarray}
where
\begin{eqnarray*}
E & = & \left\{ (k,i):\,w(E^k_i) \leq \theta\,w(Q^k_i)\right\}, \\
F & = & \left\{ (k,i):\,w(E^k_i) > \theta\,w(Q^k_i)
\mbox{ and } \sigma^k_i>\tau^k_i \right\}, \\
G & = & \left\{ (k,i):\,w(E^k_i) > \theta\,w(Q^k_i)
\mbox{ and } \sigma^k_i\leq\tau^k_i \right\},
\end{eqnarray*}
and where $\theta$ is some constant with $0<\theta<1$ which will be
chosen below.

The term $I$ is easy to estimate:
\begin{eqnarray*}
I & = & \sum_{(k,i)\in E}\alpha^{kp}\,w(E^k_i)\, \leq \, \theta
\sum_{k,i}\alpha^{kp}\,w(Q^k_i) \\
& \leq &\theta \sum_{k} \alpha^{kp}\,w\{Tf>\alpha^k\} \,
\leq \, C\,\theta \int |Tf|^p\,w\,d\mu.
\end{eqnarray*}

Let us consider the term $II$ now. By H\"older's inequality and \rf{cond2},
we obtain
\begin{eqnarray*}
II & = & \sum_{(k,i)\in F} \alpha^{kp}\,w(E^k_i) \,\leq \,
 \sum_{(k,i)\in F} w(E^k_i) \, \left(\frac{2\sigma^k_i}{\delta
 w(E^k_i)} \right)^p \\
& \leq & C\theta^{-p}\delta^{-p} \sum_{k,i} w(E^k_i)
\, \left(\frac{1}{w(Q^k_i)} \int_{U^k_i\setminus \Omega_{k+3}}
f\,\skd(w\,\chi_{E^k_i})\,d\mu\right)^p \\
& \leq & C\theta^{-p}\delta^{-p}\sum_{k,i} \frac{w(E^k_i)}{w(Q^k_i)^p}
\left(\int_{U^k_i}\! \skd(w\,\chi_{E^k_i})^{p'}\sigma\,d\mu\right)^{p/p'}
\!\!\left(\int_{U^k_i\setminus \Omega_{k+3}} \!\!\!f^pw\,d\mu\right) \\
& \leq & C\theta^{-p}\delta^{-p}\sum_{k,i}
\int_{U^k_i\setminus \Omega_{k+3}}\! f^p\,w\,d\mu.
\end{eqnarray*}
It is easy to check that the family of sets $\{U^k_i\}_i$ has finite
overlap for each $k$,
with some constant which depends on $m$. This fact implies
$$\sum_{k,i} \chi_{U^k_i\setminus\Omega_{k+3}} \leq C \sum_k \chi_{\Omega_k
\setminus\Omega_{k+3}} \leq C.$$
Therefore,
$II \leq C \int f^p\,w\,d\mu$.

Finally we have to deal with the term $III$.
Notice that $E^k_i\subset\R^d\setminus\Omega_{k+3}$ and, for $y\not\in
\Omega_{k+3}$, by Lemmas \ref{perprmax} and \ref{amesmes} we have
\begin{equation}  \label{eq7}
\sup_{x\in 2Q^{k+3}_j} s_t(y,x) \leq C\, \inf_{x\in 2Q^{k+3}_j}
\sum_{r=t-6}^{t+6} s_r(y,x)
\end{equation}
for all $t\in\Z$.
Let $H^k_i =\{j:\, Q^{k+3}_j\cap U^k_i\neq\varnothing\}$. Then,
for $j\in H^k_i$, we have
$$\sup_{x\in 2Q^{k+3}_j} \skd(w\,\chi_{E^k_i})(x) \leq
C\,\inf_{x\in 2Q^{k+3}_j}
\sum_{t=h-M-6}^{h+M+6} S_t^*(w\,\chi_{E^k_i})(x).$$
%\begin{multline*}
%\sup_{x\in 2Q^{k+3}_j} \sk(w\,\chi_{E^k_i})(x) \leq \\
%C\,\inf_{x\in 2Q^{k+3}_j}
%\sum_{t=h-M-5}^{h+M+6} S_t(w\,\chi_{E^k_i})(x).
%\end{multline*}
We set $\skk = \sum_{t=h-M-6}^{h+M+6} S_t$, and we obtain
\begin{eqnarray*}
\tau^k_i & = & \int_{U^k_i\cap \Omega_{k+3}} f\,\sk(w\,\chi_{E^k_i})\,d\mu \\
& \leq & C\, \sum_{j\in H^k_i}\, \inf_{x\in 2Q^{k+3}_j}
\skkd(w\,\chi_{E^k_i})(x)
\int_{Q^{k+3}_j}f \, d\mu \\
& \leq & C\, \sum_{j\in H^k_i} \left(\int_{2Q^{k+3}_j} \skkd(w\,\chi_{E^k_i})
\,\sigma\,d\mu \right) \left(\frac{1}{\sigma(2Q^{k+3}_j)} \int_{Q^{k+3}_j}
f\,d\mu\right).
\end{eqnarray*}
We denote $T^{k}_i =  \int_{Q^{k}_i} f\,d\mu / \sigma(2Q^{k}_i)$ and
$L^k_i = \{s:\,Q^k_s\cap U^k_i\neq \varnothing\}$. Then we have
\begin{eqnarray} \label{eq210}
\tau^k_i & \leq & C \sum_{j\in H^k_i} \left(\int_{2Q^{k+3}_j} \skkd(w\,
\chi_{E^k_i}) \,\sigma\,d\mu\right) \,T^{k+3}_j \nonumber\\
& \leq & C\sum_{s\in L^k_i}  \sum_{j:\,Q^{k+3}_j\subset Q^k_s}
\left(\int_{2Q^{k+3}_j} \skkd(w\,\chi_{E^k_i})
\,\sigma\,d\mu\right) \,T^{k+3}_j
\end{eqnarray}

We will show that
\begin{equation} \label{eq222}
\sum_{\mbox{\scriptsize $
\begin{array}{c}
(k,i)\in G \\
k\geq N_0\\
k= M_0 \mod  3  \end{array}$}}
\alpha^{kp}\,w(E^k_i) \leq C\int f^p w\,d\mu,
\end{equation}
for any $N_0$ and $M_0$.
For the rest of the proof we follow the convention that all indices
 $(k,i)$ are restricted to $k\geq N_0$ and $k= M_0 \mod  3$.

Now we will introduce {\em principal} cubes as in \cite[p.\ 540]{Sawyer} or
\cite[p.\ 804]{MW}. Let  $G_0$ be the set of indices $(k,i)$ such that $Q^k_i$
is maximal. Assuming $G_n$ already defined, $G_{n+1}$ consists of those
$(k,i)$ for which there is $(t,u)\in G_n$ with $Q^k_i \subset Q^t_u$ and
\begin{itemize}
\item[(a)] $T^k_i> 2T^t_u$,
\item[(b)] $T^r_s \leq 2T^t_u$ if $Q^k_i\subsetneq Q^r_s\subset Q^t_u$.
\end{itemize}
We denote $\Gamma = \bigcup_{n=0}^\infty G_n$, and for each $(k,i)$, we define
$P(Q^k_i)$ as the smallest cube $Q^t_u$ containing $Q^k_i$ with $(t,u)\in
\Gamma$. Then we have
\begin{itemize}
\item[(a)] $P(Q^k_i) = Q^t_u$ implies $T^k_i\leq 2T^t_u$.
\item[(b)] $Q^k_i\subsetneq Q^t_u$ and $(k,i),\, (t,u)\in \Gamma$ implies
$T^k_i>2T^t_u$.
\end{itemize}

By \rf{eq210} and the fact that $\# L^k_i\leq C$, we get
\begin{eqnarray*}
\lefteqn{\sum_{(k,i)\in G} \alpha^{kp}\,w(E^k_i) \leq
\sum_{(k,i)\in G} w(E^k_i)\left( \frac{2\tau^k_i}{\delta w(E^k_i)}\right)^p
 \leq \, C\,\sum_{k,i}
\frac{w(E^k_i)}{w(Q^k_i)^p}\, (\tau^k_i)^p} &&  \\
& \leq & \mbox{}\!\!\! C \sum_{k,i} \sum_{s\in L^k_i}
\frac{w(E^k_i)}{w(Q^k_i)^p}\,
\Biggl[ \!\!\sum_{\mbox{\scriptsize $
\begin{array}{c} j: (k+3,j)\not\in\Gamma\\
Q^{k+3}_j\subset Q^k_s\\
P(Q^{k+3}_j) = P(Q^k_s)
\end{array}$}} \!\!\!\!\!\!
\left(\int_{2Q^{k+3}_j} \skkd(w\,\chi_{E^k_i})
\,\sigma\,d\mu\right) \,T^{k+3}_j \Biggr]^p \\
&& \mbox{} + C\,\sum_{k,i} \frac{w(E^k_i)}{w(Q^k_i)^p}
\Biggl[\sum_{j\in H^k_i:\,(k+3,j)\in\Gamma}
\left(\int_{2Q^{k+3}_j} \!\!\skkd(w\,\chi_{E^k_i})
 \,\sigma\,d\mu\right) \,T^{k+3}_j \Biggr]^p \\
& = &  IV + V.
\end{eqnarray*}

Let us estimate the term $IV$ first. Notice that if
$(k+3,j)\not\in\Gamma$, then $Q^{k+3}_j\neq P(Q^{k+3}_j)$. As a consequence,
$\ell(Q^{k+3}_j)\leq \ell(P(Q^{k+3}_j))/2$, and $2Q^{k+3}_j\subset
\frac{4}{3}P(Q^{k+3}_j)$. Taking into account this fact, the finite
overlap of the cubes $Q^{k+3}_j$ (for a fixed $k$), and \rf{cond1},
for any $(t,u)\in\Gamma$ we get
\begin{eqnarray*}
\lefteqn{
\sum_{k,i} \sum_{s\in L^k_i:\, P(Q^k_s)= Q^t_u} \!
\frac{w(E^k_i)}{w(Q^k_i)^p}
\Biggl[\!\! \sum_{\mbox{\scriptsize $
\begin{array}{c} j: (k+3,j)\not\in\Gamma\\
Q^{k+3}_j\subset Q^k_s\\
P(Q^{k+3}_j) = Q^t_u
\end{array}$}}  \!\!\!\!\!
\left(\int_{2Q^{k+3}_j} \!\!\skkd(w\,\chi_{E^k_i})
\,\sigma\,d\mu\right) \,T^{k+3}_j \Biggr]^p} && \\
& \leq & C\sum_{k,i}  \sum_{s\in L^k_i:\, P(Q^k_s)= Q^t_u} \!\!\!w(E^k_i)
\left(\frac{1}{w(Q^k_i)} \int_{\frac{4}{3}Q^t_u}
\skkd(w\,\chi_{Q^k_i}) \,\sigma\,d\mu\right)^p (2T^t_u)^p \\
& \leq & C\,(T^t_u)^p \sum_{k,i}  \sum_{s\in L^k_i:\, P(Q^k_s)= Q^t_u}\!\!
w(E^k_i) \left(\frac{1}{w(Q^k_i)} \int_{Q^k_i}
\skk(\sigma\,\chi_{\frac{4}{3}Q^t_u}) \,w\,d\mu\right)^p \\
& \leq &  C\,(T^t_u)^p \int M^d_{w}(N(\sigma\,\chi_{Q^t_u}))^p\,w\,d\mu \\
& \leq & C\,(T^t_u)^p \int
N(\sigma\,\chi_{\frac{4}{3}Q^t_u})^p\,w\,d\mu
\,\leq \,C\,(T^t_u)^p\,\sigma(2Q^t_u),
\end{eqnarray*}
where we have denoted by $M^d_{w}$ the dyadic maximal operator with
respect to $w$. Thus,
\begin{equation} \label{eqiv}
IV \leq C\sum_{(t,u)\in\Gamma} \sigma(2Q^t_u) (T^t_u)^p.
\end{equation}

Let us estimate the term $V$. By H\"older's inequality and \rf{cond2},
for a fixed $(k,i)$,
\begin{eqnarray*}
\lefteqn{\frac{w(E^k_i)}{w(Q^k_i)^p}
\Biggl[\sum_{j\in H^k_i:\,(k+3,j)\in\Gamma}
\left(\int_{2Q^{k+3}_j} \!\!\skkd(w\,\chi_{E^k_i})
\,\sigma\,d\mu\right) \,T^{k+3}_j \Biggr]^p} && \\
& \leq & \frac{w(E^k_i)}{w(Q^k_i)^p}
\Biggl[\sum_{j\in H^k_i} \sigma(2Q^{k+3}_j)^{-p'/p}
\left(\int_{2Q^{k+3}_j} \skkd(w\,\chi_{E^k_i})
\,\sigma\,d\mu\right)^{p'}\Biggr]^{p/p'} \\
&&\mbox{} \times \Biggl[\sum_{j\in H^k_i:\,(k+3,j)\in\Gamma}
\sigma(2Q^{k+3}_j)\, (T^{k+3}_j)^p\Biggr] \\
& \leq & \frac{w(E^k_i)}{w(Q^k_i)^p}
\Biggl[\sum_{j\in H^k_i}
\int_{2Q^{k+3}_j} \skkd(w\,\chi_{E^k_i})^{p'}
\,\sigma\,d\mu\Biggr]^{p/p'} \\
&&\mbox{} \times
\,\Biggl[\sum_{j\in H^k_i:\,(k+3,j)\in\Gamma}
\sigma(2Q^{k+3}_j)\, (T^{k+3}_j)^p\Biggr] \\
& \leq & C\sum_{j\in H^k_i:\,(k+3,j)\in\Gamma}
\sigma(2Q^{k+3}_j)\, (T^{k+3}_j)^p.
\end{eqnarray*}
Summing over $(k,i)$, since any cube $Q^{k+3}_j$ occurs at most $C$ times
in the resulting sum, we get
\begin{equation}  \label{eqv}
V \leq C\,\sum_{(t,u)\in\Gamma} \sigma(2Q^t_u) (T^t_u)^p.
\end{equation}

Notice that for each $(t,u)$ we can write
\begin{eqnarray*}
\sigma(2Q^t_u) (T^t_u)^p & = &\sigma(Q^t_u) (T^t_u)^{p-1}
\frac{1}{\sigma(Q^t_u)} \int_{Q^t_u} f\sigma^{-1}\,\sigma\,d\mu \\
& =: & \sigma(Q^t_u) (T^t_u)^{p-1} m_{\sigma,Q^t_u}(f\sigma^{-1}).
\end{eqnarray*}
We have obtained
\begin{eqnarray*}
IV  + V & \leq & C\,\sum_{(t,u)\in\Gamma} \sigma(Q^t_u) (T^t_u)^{p-1}
m_{\sigma,Q^t_u}(f\sigma^{-1}) \\
& = & C\,\int \Biggl(\sum_{(t,u)\in\Gamma} (T^t_u)^{p-1}
m_{\sigma,Q^t_u}(f\sigma^{-1}) \,\chi_{Q^t_u}(x)
\Biggr)\,\sigma(x)\,d\mu(x).
\end{eqnarray*}
Notice that for any fixed $x$ we have
\begin{eqnarray*}
\sum_{(t,u)\in\Gamma} (T^t_u)^{p-1}m_{\sigma,Q^t_u}(f\sigma^{-1})
 \,\chi_{Q^t_u}(x) & \leq & C\,\sup_{(t,u)
\in\Gamma:\,x\in Q^t_u} (T^t_u)^{p-1}\, M^d_\sigma(f\sigma^{-1})(x)\\
& \leq & C\, M^d_\sigma(f\sigma^{-1})(x)^p.
\end{eqnarray*}
Therefore,
$$IV  + V \leq  C\,\int  M^d_\sigma(f\sigma^{-1})^p\,\sigma\,d\mu \leq
 C\,\int (f\sigma^{-1})^p\,\sigma\,d\mu = C\,\int f^p\,w\,d\mu,$$
which yields \rf{eq222}.
Thus, by \rf{eq123},
$$\int |Tf|^p\,w\,d\mu \leq C\,(I+II+III) \leq
C\theta\int |Tf|^p\,w\,d\mu +  C\,\int f^p\,w\,d\mu.$$
We only have to choose $\theta$ small enough, and we are done.
\end{proof}

% ****************************************************************************
% ****************************************************************************

To complete the proof of the implication (d) $\Rightarrow$ (c) of
Lemma \ref{main}, it remains to prove the following result.

\begin{lemma} \label{lemaclau}
With the notation and assumptions of Lemma \ref{dif2},
\rf{fonam} holds. That is,
$w\biggl( \bigcup_{j\geq k} B_j \biggr) \leq \eta\, w(\Omega_k)$, with
$0<\eta<1$
\end{lemma}

Before proving the lemma, a remark:

\begin{rem} \label{besico}
Besicovitch's Covering Theorem asserts that if $A\subset\R^d$ is
bounded and there exists some family of cubes ${\mathcal Q} =
\{Q_x\}_{x\in A}$, with each $Q_x$ centered at $x$, then there exists
some finite or countable family of cubes $\{Q_{x_i}\}_i\subset
{\mathcal Q}$ which covers $A$ with finite overlap. That is,
$\chi_A\leq \sum_i \chi_{Q_{x_i}} \leq C$, with $C$ depending only on
$d$.

We are going to show that the covering $\{Q_{x_i}\}_i$ can be chosen
so that the following property holds too:
\begin{equation} \label{besic2}
\mbox{\em If $z\in A\cap Q_{x_i}$
for some $i$, then $\ell(Q_z) \leq 4 \ell(Q_{x_i})$.}
\end{equation}
Indeed, for each $x\in A$, let $R_x$ be some cube of the type $Q_y$,
$y\in A$, with $x\in\frac{1}{2}R_x$ and such that
$$\ell(R_x) \geq
\frac{99}{100}\,\sup_{y:x\in\frac{1}{2}Q_y}\ell(Q_y).$$
Now we will apply Besicovitch's Covering Theorem to the family of
cubes $\{R_x\}_{x\in A}$. Let us remark that the Theorem of
Besicovitch also holds for the family $\{R_x\}_{x\in A}$ because,
although the cubes $R_x$ are not centered at $x$, we still have
$x\in\frac{1}{2}R_x$ (see \cite{Morse} or \cite[p.\ 6-7]{Guzman}, for example).
So there exists some
finite or countable family $\{R_{x_i}\}_i$ which covers $A$ with
finite overlap. Notice that $\{R_{x_i}\}_i\subset {\mathcal Q}$, and
if $z\in R_{x_i}\cap A$, then $\ell(Q_z)\leq
4\ell(R_{x_i})$. Otherwise, $x_i\in \frac{1}{2}Q_z$ and $\ell(Q_z)>
4\ell(R_{x_i})$, which contradicts the definition of $R_{x_i}$.

It is worth comparing this version of Besicovitch Covering Theorem
with the version of Wiener's Covering Lemma \ref{vitali}. Notice that the
statement (3) of Lemma \ref{vitali} and \rf{besic2} look quite similar.
\end{rem}

% ***************************************************************************

\vspace{2mm}
\begin{proof}[{\bf Proof of Lemma \ref{lemaclau}}]
We use the same notation as in the proof of the preceding lemma.

Let $x\in B_j$ and take $Q^j_i$ containing $x$
(recall $B_j\subset\Omega_{j+2}\subset\Omega_j)$, with $Q^j_i\in \AD_{h}$.
By \rf{andso}, we have $N(f\,\chi_{U^j_i})(x)\geq \ve_0 \,\alpha^j$
for some $\ve_0>0$ depending on $\alpha,\,\beta,\,m$.
It is easily seen that this implies that $S_t(f\,\chi_{U^j_i})(x)\geq
\ve_0 \,\alpha^j$ for $t\geq h-M$, where $M$ is some positive constant which
depends on $\ell(U^j_i)/\ell(Q^j_i)$. Also, by the definition of $B_j$,
$S_t (f\,\chi_{U^j_i})(x)\leq \delta\alpha^j$ for $h-M\leq t\leq h+M$.
So, if we choose $\delta\leq \ve_0$ and $M\geq10$, then
\begin{equation} \label{vv0.5}
\sup_{t\geq h+10} S_t f(x)>\ve_0\,\alpha^j
\end{equation}

We denote $A_k := \bigcup_{j\geq k} B_j$. For a fixed $x\in A_k$, let $r$
be the least integer such that $r\geq k$ and $x\in B_r$. There exists some cube
$Q^r_i$ containing $x$, with $Q_i^r\in\AD_{h}$ for some $h$. Since
$S_{h+5}(f\,\chi_{U^r_i})(x) \leq \delta\alpha^r$, by Lemma \ref{guai}
there exists some doubling cube $P_x\in\AD_{h+5,h+4}$ centered at $x$ such
that
\begin{equation} \label{vv1}
\frac{1}{\mu(2P_x)} \int_{2P_x} f\,d\mu \leq C\delta\alpha^r.
\end{equation}
Now, by Besicovitch's Covering Theorem, we can find some family of cubes
$\{P_{x_s}\}_{s}$ (with $x_s\in A_k$) which covers $A_k$ with finite
overlap. Moreover, we assume that the covering has been chosen so that the
property \rf{besic2} holds.

Given any $\rho$ with $0<\rho<1$, we will show that if $\delta$ is small
enough, then
\begin{equation} \label{vv2.5}
\mu(A_k\cap P_{x_s}) \leq \rho\,\mu(P_{x_s})
\end{equation}
for all $s$.

Let $P_{x_0}$ some fixed cube from the family $\{P_{x_s}\}_s$, and let $r_0$
be the least integer such that $x_0\in B_{r_0}$. First we will see that
\begin{equation} \label{vv2.55}
\mu\biggl(\bigcup_{j\geq r_0} B_j\cap P_{x_0}\biggr) \leq
\frac{\rho}{2}\,\mu(P_{x_0}).
\end{equation}
If $z\in B_j\cap P_{x_0}$ for some $j\geq r_0$ and $z\in Q^j_i$, then by
\rf{vv0.5} we have
\begin{equation}  \label{vv1.99}
\sup_{j\geq h+10} S_jf(z) > \ve_0\,\alpha^j,
\end{equation}
where $h$ is so that $Q^j_i\in\AD_{h}$.
Let us denote by $Q^{r_0}_{i_0}$ the Whitney cube of $\Omega_{r_0}$ containing
$x_0$, with $Q^{r_0}_{i_0}\in\AD_{h_0}$. Since
$\Omega_j\subset\Omega_{r_0}$, we have $\ell(Q^j_i)\leq
C_{26}\,\ell(Q^{r_0}_{i_0})$, and so $h\geq h_0-1$. In fact, if
$C_{26}$, which depends on $d$,
is very big, then we should write $h\geq h_0-q$, where $q$ is
some positive integer big enough, depending on $C_{26}$.
The details of the required modifications in this case are left to the
reader. From \rf{vv1.99}, we get
\begin{equation}  \label{vv2}
\sup_{j\geq h_0+9} S_jf(z) > \ve_0\,\alpha^j\geq\ve_0\,\alpha^{r_0}.
\end{equation}
For $j\geq h_0+9$ and $z\in P_{x_0}$, we have $\supp(s_j(z,\cdot))\subset
2P_{x_0}$, because $P_{x_0}\in \AD_{h_0+5,h_0+4}$. Thus \rf{vv2} implies
$N(f\,\chi_{2P_{x_0}})(z) > \ve_0\,\alpha^{r_0}$,
and then, from the weak $(1,1)$ boundedness of $N$, by \rf{vv1}, and
because $P_{x_0}$ is doubling, we obtain
\begin{eqnarray} \label{vv5}
\mu\biggl(\bigcup_{j\geq r_0} B_j \cap P_{x_0}\biggr) & \leq &
\mu\left\{z\in P_{x_0}:\, N(f\,\chi_{2P_{x_0}})(z) >
  \ve_0\,\alpha^{r_0}\right\}
\nonumber\\
& \leq & \frac{C}{\ve_0\,\alpha^{r_0}} \int_{2P_{x_0}}f\, d\mu
% \leq  \frac{C}{\ve_0\,\alpha^{r_0}}\, \delta\,
%\alpha^r\,\mu(2P_{x_0})
\,\leq \, C\,\ve_0^{-1}\,\delta\,\mu(P_{x_0}).
\end{eqnarray}
So \rf{vv2.55} holds if $\delta$ is sufficiently small.

Now we have to estimate $\mu\Bigl(\bigcup_{k\leq j\leq r_0-1} B_j\cap
P_{x_0}\Bigr)$. If $z\in B_j\cap P_{x_0}$, then
$\ell(P_z)\leq 4\,\ell(P_{x_0})$, by \rf{besic2}.
Recall also that
$P_{x_0}\in\AD_{h_0+5,h_0+4}$. As a consequence, we deduce
$P_z\in \AD_{+\infty,h_0+3}$. Moreover, we have
$P_{x_0}\subset\{Tf>\alpha^{r_0}\}$, and
so $Nf(z) >C_{27}\alpha^{r_0}$, with $C_{27}>0$. Since by \rf{vv****} we have
$$N(f\,\chi_{\R^d\setminus U^j_i})(z) \leq C_{28}\,\alpha^{j},$$
we obtain
$$N(f\,\chi_{U^j_i})(z) > C_{27}\alpha^{r_0} -  C_{28}\,\alpha^{j}
\geq  \frac{C_{27}}{2}\,\alpha^{r_0},$$
assuming $j\leq r_0-r_1$, where $r_1$ is some positive integer which depends on
$C_{27}$ and $C_{28}$. Recall also that the fact that $z\in B_j$ yields
\begin{equation} \label{vv4}
S_t(f\,\chi_{U^j_i})(z) \leq \delta\,\alpha^j\qquad \mbox{for $h_1-10
  \leq t\leq h_1+10$},
\end{equation}
where $h_1$ is given by $Q^j_i\in\AD_{h_1}$. If we choose
$\delta$ small enough, then $\delta\,\alpha^j\leq C_{27}\,\alpha^{r_0}/2$
and, for $j\leq r_0-r_1$, \rf{vv4} implies
\begin{equation}  \label{vv*1}
S_t(f\,\chi_{U^j_i})(z) > \frac{C_{27}}{2}\,\alpha^{r_0} \qquad\mbox{for some
$t\geq h_1+10$.}
\end{equation}

On the other hand, if $r_0-r_1<j<r_0$, then by \rf{vv0.5} we have
\begin{equation}  \label{vv*2}
S_tf(z) > \ve_0\,\alpha^{r_0-r_1} \geq C_{29}\,\alpha^{r_0}
\qquad\mbox{for some $t\geq h_1+10$,}
\end{equation}
with $C_{29}>0$.

In any case, from the fact that $P_z\in\AD_{h_1+5,h_1+4}$ we
deduce $h_1\geq h_0-2$, and so
$\supp(s_t(z,\cdot))\subset 2P_{x_0}$ for $t\geq h_1+10$. Thus, from
\rf{vv*1} and \rf{vv*2} we get
$$N(f\,\chi_{2P_{x_0}})(z) \geq \min(C_{27}/{2},\, C_{29})\,\alpha^{r_0}$$
for any $j$ with $k\leq j <r_0$. If we take $\delta$ small enough,
operating as in \rf{vv5}, we obtain
$$
\mu\biggl(\bigcup_{k\leq j< r_0} B_j \cap P_{x_0}\biggr)
\leq C\delta \mu(P_{x_0}) \leq \frac{\rho}{2}\, \mu(P_{x_0}),
$$
which together with \rf{vv2.55} implies \rf{vv2.5}.

By \rf{vv2.5} and Lemma \ref{cubspetits},
using the $Z_\infty$ condition for $w$, we get
$w(2Q^k_i\setminus A_k) \geq \tau\,w(Q^k_i)$ for each Whitney cube
$Q^k_i\in\Omega_k$.
By the finite overlap of the cubes $2Q^k_i$, we obtain
$$\tau\,w(\Omega_k) \leq \tau \sum_i w(Q^k_i) \leq \sum_i
w(2Q^k_i\setminus A_k) \leq C_{30} \,w(\Omega_k\setminus A_k).$$
Therefore,
$$w(A_k) \leq (1-C_{30}^{-1}\tau)\,w(\Omega_k) =: \eta\,w(\Omega_k).$$
\end{proof}

% ***************************************************************************
% ***************************************************************************
% ***************************************************************************

\section{The general case}

In this section we consider the case where not all the cubes
$Q_{x,k}\in\DD$ are transit cubes.

If $\R^d$ is an initial cube but there are no stopping cubes, then
the arguments in Sections 5--10 with some minor
modifications are still valid.

If there exist stopping cubes, some problems arise because the functions
$S_k\chi_{\R^d}$ are not bounded away from zero, in general.
As a consequence, the property $Z_\infty$ has to be modified.
Indeed, notice that if we set $A:=\R^d$ and $Q$ is some cube which contains
stopping points, then \rf{hsinfinit} may fail, and so the $Z_\infty$
condition is useless in this case.

The new formulation of the $Z_\infty$ property is the following.
For $k\in\Z$, we denote $\ST_k:=
\{x\in\supp(\mu): Q_{x,k}\mbox{ is a stopping cube}\}$. Notice by the way
that $S_jf(x)=0$ for $j\geq k+2$ and $x\in\ST_k$.

\begin{definition}  \label{sinfinit'}
We say that $w$ satisfies the $Z_\infty$ property if there exists some
constant
$\tau>0$ such that for any cube $Q\in \AD_{k}$ and any set
$A\subset \R^d$ with $Q\cap \ST_{k+3}\subset A$, if
\begin{equation}  \label{hsinfinit'}
S_{k+3}\chi_A(x)\geq 1/4\qquad \mbox{for all $x\in Q\setminus\ST_{k+3}$,}
\end{equation}
then $w(A\cap 2Q) \geq \tau\, w(Q)$.
\end{definition}

With this new definition, Lemma \ref{lsinfinit} still holds. The new proof is
a variation of the former one. On the other hand, Lemma
\ref{cubspetits} changes. Let us state the new version:

\begin{lemma}  \label{cubspetits'}
Suppose that $w$ satisfies the $Z_\infty$ property.
Let $Q\in\AD_{h}$ and $A\subset\R^d$ be such that $A\cap
Q\cap\ST_{h+4}=\varnothing$. Let $\{P_i\}_i$ be
a family of cubes with finite overlap such that
$A\cap\frac{3}{2}Q\subset
\bigcup_i P_i$, with $P_i\in \AD_{+\infty,h+4}$ and $\ell(P_i)>0$
for all $i$. There exists some constant $\delta>0$ such that if
$\mu(A \cap P_i) \leq \delta \,\mu(P_i)$ for each $i$, then
\begin{equation} \label{bm0.5'}
w(2Q\setminus A) \geq \tau\,w(Q),
\end{equation}
for some constant $\tau>0$ (depending on $Z_\infty$).
If, moreover, $w(2Q) \leq C_{11}\,w(Q)$, then
\begin{equation} \label{bm1'}
w(A\cap 2Q) \leq (1-C_{11}^{-1}\tau)\,w(2Q).
\end{equation}
\end{lemma}

The proof is analogous to the proof of Lemma \ref{cubspetits}, and it is left
for the reader.

The results stated in the other lemmas in Sections 5--10 remain true
in the new situation.
However, the use of the $Z_\infty$ condition is basic in the proofs of Lemma
\ref{lemagli}, the implication (e) $\Rightarrow$ (c) of Lemma \ref{mainw},
Lemma \ref{fonam'}, Lemma \ref{dif2}, and Lemma \ref{lemaclau}. Below we will
describe the changes required in the arguments.
In the rest of the lemmas and results, the proofs and arguments either
are identical or require only some minor modifications (which are left for
the reader again).

% ***************************************************************************

\vspace{3mm}
\noindent{\em - Changes in the proof of Lemma \ref{lemagli}.}
The proof is the same until \rf{nn*1}, which still holds.
Given $Q_i\in\AD_k$, it is
easily seen that if $y\in \ST_{k+3}\cap Q_i$, then $T_*(f\chi_{3Q_i})(y)
\leq C_{31} Nf(y)$. By \rf{nn*1}, if we choose $\delta<
\ve/2C_{31}$, then $A_\lambda\cap Q_i\cap\ST_{k+3}=\varnothing$.

On the other hand, now the estimate \rf{claim1} is valid for $y\in Q_i
\setminus \ST_{k+3}$. Then we deduce $S_{k+3}\chi_{2Q_i\setminus A_\lambda}
(y)>\frac14$ for $y\in Q_i\setminus\ST_{k+3}$, and by the $Z_\infty$,
condition we get $w(2Q_i\setminus A_\lambda)\geq\tau w(Q_i)$. Arguing as in
\rf{nn*2}, we obtain $w(A_\lambda) \leq \rho,w(\Omega_\lambda)$.

% ***************************************************************************

\vspace{3mm}
\noindent{\em - Changes in the proof of
the implication (e) $\Rightarrow$ (c) of Lemma \ref{mainw}.}
The sets $\Omega_\lambda,\,G_\lambda$ and $B_\lambda$ are defined
in the same way. The estimates for $w(Q_i\cap G_\lambda)$ are the same.

As shown in \rf{nn*10}, if $z\in B_\lambda\cap Q_i$, with $Q_i\in\AD_h$, then
$S_k(f\chi_{U_i})(z)\geq C_{16}\lambda\neq0$ for some $k\geq h+6$.
This implies $z\not\in \ST_{h+4}$.
Now the arguments used to prove that
$w(B_\lambda)\leq \eta_1 w(\Omega_\lambda)$
are still valid, because $B_\lambda\cap Q_i\cap \ST_{h+4}=\varnothing$.

% ***************************************************************************

\vspace{3mm}
\noindent{\em - Changes in the proof of Lemma \ref{fonam'}.}

\vspace{1mm}
\noindent{\em The construction.} The construction is basically the same. The
only difference is that now we must be careful because the cubes $Q^1_x$ (and
$Q_{x}^k$ for $k>1$) may fail to exist due to the existence of stopping points.
In the first step of the construction ($k=1$), we circumvent this problem as
follows. If $x\in R^1_j\setminus \ST_{g^1_j+18}$,
then we take a $\mu$-$\sigma$-$(100,\beta)$-doubling cube $Q^1_x\in
\AD_{g^1_j+16}$. If $x\in R^1_j\cap\ST_{g^1_j+18}$, we write $x\in\AS_1$.
We consider a Besicovitch covering of $Q_0\cap\Omega_0\setminus\AS_1$ with
this type of cubes: $Q_0\cap\Omega_0\setminus\AS_1\subset \bigcup_{i\in I_1}
Q^1_i$, and we set $A_1:=\bigcup_{i\in I_1} Q^1_i$.
We operate in an analogous way at each step $k$ of the construction.

\vspace{1mm}
\noindent{\em The estimate of $\int_{Q_0} |N^h\sigma|^p w\,d\mu$.}
Here there are little changes too. Equation \rf{mmm10} is proved inductively
in the same way. Let us see the required modifications in the first step.
The definition of $B_0$ is different now: $B_0:=\{x\in Q_0\setminus\ST_{h+5}:
S_{h+3}\sigma(x)\leq\ve\lambda_0\}$. With this new definition,
\rf{mmm*4} holds. On the other hand, notice that $\int_{Q_0
\cap\ST_{h+5}}|N^{h+20}\sigma|^p w\,d\mu=0$, since $N^{h+20}\sigma(x)=0$ if
$x\in\ST_{h+5}$.

The definition of $D_0$ does not change, and all the other estimates remain
valid. In particular, \rf{nn*21} holds now too, because $N^{H^1_x+20}
\sigma(x)=0$ if $x\in\AS_1$ (recall that the definition of $A_1$ has changed).

The changes required at each step $k$ are analogous.

\vspace{1mm}
\noindent{\em The estimate of $\sum_k \sigma(A_k)$.} The former arguments
remain valid in the new situation.

% ***************************************************************************

\vspace{3mm}
\noindent{\em - Changes in the proof of Lemma \ref{dif2} and Lemma
\ref{lemaclau}.}
The proof of Lemma \ref{dif2} does not change.
In the arguments for Lemma \ref{lemaclau}, we have to take into account that
if $x\in B_k$ and $\delta$ is small enough, then $x\not\in\ST_{h+M-1}$. Indeed,
if $x\in Q_i$, with $Q_i\in\AD_h$, then we have
$T(f\chi_{U_i^k})(x)\geq\frac{\alpha-1}\alpha \alpha^{k+2}$, and
$S_j(f\chi_{U_i^k})(x)\leq \delta\alpha^k$ for $h-M\leq j\leq h+M$.
These inequalities imply $S_j(f\chi_{U_i^k})(x) >\ve_0\alpha^k\neq0$ for
some $j\geq M+1$ if $\delta$ is small enough.
In particular, $x\not\in \ST_{h+M-1}$.

If we assume $M\geq20$, for instance, then all the cubes $P_x$ that appear
in the proof of Lemma \ref{lemaclau} exist and are transit cubes, and
the same estimates hold.

% ***************************************************************************
% ***************************************************************************
% ***************************************************************************

\section{Relationship with $\rbmo(\mu)$ and final remarks}

Let us recall one of the equivalent definitions of the space $\rbmo(\mu)$
introduced in \cite{Tolsa2}.
We say that $f\in L^1_{loc}(\mu)$ belongs to $\rbmo(\mu)$ if there
exists a collection of numbers $\{f_Q\}_{Q\subset\R^d}\subset \R$
such that
$$\int_Q |f(x)-f_Q|\,d\mu(x) \leq C_f\,\mu(2Q)$$
for each cube $Q\subset \R^d$ and
\begin{equation} \label{defrbmo2}
|f_Q - f_R|\leq C_f (1+\delta(Q,R))
\end{equation}
for all the cubes $Q,R$ with $Q\subset R$.
The optimal constant $C_f$ is the $\rbmo(\mu)$ norm of $f$, which we
denote by $\|f\|_*$.

Let $1<p<\infty$. In general, if $w\in Z_p$, then $\log w\not\in Z_p$.
This follows easily from Example \ref{ex3}. Indeed, in this case
it can be checked that $\delta(0,I_k) \leq C$ uniformly on $k$.
As a consequence, for all $f\in\rbmo(\mu)$, the numbers $f_{I_k}$ are bounded
uniformly on $k$. If moreover $f$ is constant on each interval $I_k$, then we
deduce $f\in L^\infty(\mu)$. However, the weight $w_0$ of Example \ref{ex3}
is constant on each interval $I_k$ and it is not a bounded function, and
the same happens with $\log w_0$.

On the other hand, if $f\in\rbmo(\mu)$, then there exists some
$\ve>0$ depending on $\|f\|_*,\,p$ such that $e^{\ve f}\in Z_p$. To prove
this,
first we will show in the following proposition that a weight of the type
$e^{\ve f}$, with $f\in\rbmo(\mu)$,
satisfies a (rather strong) property in the spirit of the classical
$A_p$ condition.

\begin{propo}
Let $1<p<\infty$.
If $f\in\rbmo(\mu)$ and $\ve=\ve(\|f\|_*,p)>0$ is small enough, then
\begin{equation} \label{AP}
\frac{1}{\mu(2Q)} \int_Q e^{\ve f}d\mu \,\cdot\, \left[ \frac{1}{\mu(2R)}
\int_R e^{-\ve f p'/p} d\mu\right]^{p/p'} \leq C_{32}\,
e^{C_{33}\delta(Q,R)},
\end{equation}
for all the cubes $Q,R$ with $Q\subset R$ or $R\subset Q$, where
$C_{32},C_{33}$ are positive constants depending on $n,d,C_0$.
\end{propo}

\begin{proof}
The funcions $f\in\rbmo(\mu)$ satisfy an inequality of John-Nirenberg type
(see \cite[Theorem 3.1]{Tolsa2}), which implies that for some constants
$C_{34}, C_{35}$ and any cube $Q$ and $\lambda>0$ we have
$$\int_Q \exp(C_{34}|f(x)-f_Q|/\|f\|_*)\, d\mu(x) \leq C_{35}\,\mu(2Q).$$
If we take $\ve\leq C_{34}\min(1,p/p')/\|f\|_*)$ and we use \rf{defrbmo2},
we deduce \rf{AP}.
\end{proof}

\begin{rem}
For $1<p<\infty$, in Lemma \ref{mainw}, the statement (e) can be substituted
by the following weaker assumption:
\begin{enumerate}
\item[(e')] For all $k\in\Z$ and all cubes $Q$,
\begin{equation} \label{cond1''}
\int_Q |S_k(w\chi_Q)|^{p'}\, \sigma\,d\mu \leq C\,w(2Q),
\end{equation}
with $C$ independent of $k$ and $Q$.
\end{enumerate}
To see this, only some minor changes (which are left for the reader) in the
proof of Lemma \ref{mainw} are required.

Since (e) and (e') in Lemma \ref{mainw} are equivalent, we deduce that the
statement (e) of Lemma \ref{main} can be weakened in the analogous way:
We only need to compute both integrals over $Q$, and on the right hand side
$Q$ can be substituted by $2Q$.
\end{rem}

\begin{theorem} \label{rbmo}
Let $1<p<\infty$.
If $f\in\rbmo(\mu)$ and $\ve= \ve(\|f\|_*,p)>0$ is small enough,
then $e^{\ve f}\in Z_p$.
\end{theorem}

\begin{proof}
By the preceding remark, we only need to show that $w:=e^{\ve f}$
satisfies \rf{cond1''} and its corresponding dual estimate.
Moreover, for simplicity we will assume that there are no stopping cubes.

Let us see that \rf{cond1''} holds any given cube $Q\in\AD_h$. We may assume
$k\geq h-3$, since
$S_k(w\chi_Q) (x) \leq C\sum_{i=-3}^3 S_{h+i}(w\chi_Q) (x)$
for $x\in Q\cap\supp(\mu)$.

For each $x\in Q\cap\supp(\mu)$, let $R_x$ be a doubling cube centered at $x$,
with $R_x\in\AD_{k+10,k+9}$. Let $\bigcup_i R_i\supset Q\cap \supp(\mu)$ be a
Besicovitch covering of $Q\cap \supp(\mu)$ with this type of cubes. Notice that
$R_i\subset 2Q$ for all $i$. Let $Q_x\in\AD_{k,k-2}$ be some cube centered
at $x$. If $x\in R_i$, then $R_i\subset Q_x$ because $\ell(R_i)\ll \ell(Q_x)$.
Since $\delta(R_i,Q_x)\leq C$, by \rf{AP} we have
$$\frac{1}{\mu(2Q_x)} \int_{Q_x} w\,d\mu \,\cdot\, \left[ \frac{1}{\mu(R_i)}
\int_{R_i} \sigma\, d\mu\right]^{p/p'} \leq C.$$
Taking a suitable mean over cubes $Q_x$ centered at $x$ (as in the proof of
Lemma \ref{guai}), we obtain
$$S_k w(x)\cdot \bigl[m_{R_i}
(\sigma)\bigr]^{p/p'} \leq C,$$
for all $x\in Q$. Then we get
\begin{eqnarray*}
\int_Q |S_k(w\chi_Q)|^{p'}\, \sigma\,d\mu & \leq &  \sum_i
\int_{R_i} |S_k(w)|^{p'}\, \sigma\,d\mu \\
& \leq & C\sum_i \frac{\sigma(R_i)}{(m_{R_i}\sigma)^p}
\,= \,C\sum_i \frac{\mu(R_i)}{(m_{R_i}\sigma)^{p-1}}.
\end{eqnarray*}
By H\"older's inequality, $1\leq m_{R_i}w \cdot (m_{R_i}\sigma)^{p-1}$.
Thus,
$$\int_Q |S_k(w\chi_Q)|^{p'}\, \sigma\,d\mu \leq C\sum_i m_{R_i}w\cdot
\mu(R_i) = C \sum_i w(R_i) \leq C\,w(2Q).$$
The estimate dual to \rf{cond1''} is proved in an analogous way.
\end{proof}

\vspace{3mm}
We will finish with some remarks and open questions:
\begin{rem}

\begin{itemize}

\item[(a)]
Using Lemma \ref{fonam'} and modifying a little the proof
  of the implication (e) $\Rightarrow$ (c) of Lemma \ref{mainw} one
  can show that $w\in Z_p^{weak}$ if and only if there exists some
  $\lambda>1$ such that
$\int N(w\chi_Q)^p\, \sigma\,d\mu \leq C\,w(\lambda Q)$
for all cubes $Q$. We don't know if this holds with $\lambda=1$ too.

\item[(b)] We don't know if $Z_p = Z_p^{weak}$.

\item[(c)] In the case $p=1$, an statement such as (e) in Lemma \ref{mainw}
is missing. We don't know if there is a reasonable substitute.
\end{itemize}
\end{rem}

% ***************************************************************************
% ***************************************************************************
% ***************************************************************************

% ***************************************************************************
% ***************************************************************************
% ***************************************************************************

% ***************************************************************************
% ***************************************************************************
% ***************************************************************************

% ***************************************************************************
% ***************************************************************************

\end{document}